\documentclass[twoside,reqno]{birkart}

\usepackage{amssymb}
\usepackage{enumerate}

\input epsf

\newtheorem{theorem}{Theorem}[section]
\newtheorem{lemma}[theorem]{Lemma}
\newtheorem{corollary}[theorem]{Corollary}
\newtheorem{proposition}[theorem]{Proposition}

\newtheorem{remark}[theorem]{Remark}
\newtheorem{examples}[theorem]{Examples}
\newtheorem{example}[theorem]{Example}
\newtheorem{conjecture}[theorem]{Conjecture}
\newtheorem{definition}[theorem]{Definition}
\newtheorem{problem}[theorem]{Problem}

\newcommand{\Z}{{\mathbb Z}}
\newcommand{\R}{{\mathbb R}}
\newcommand{\Q}{{\mathbb Q}}
\newcommand{\C}{{\mathbb C}}
\newcommand{\lra}{\longrightarrow}
\newcommand{\lms}{\longmapsto}
\newcommand{\id}{\mathop\mathrm{Id}}

\newcommand{\Ndash}{\nobreakdash--}
\newcommand{\guio}[1]{\nobreakdash-\hspace{0pt}#1}

\begin{document}

\title[Multiple $\zeta$-Values,  Galois Groups
and Geometry of Modular Varieties]
{Multiple $\zeta$-Values,  Galois Groups,
and \\ Geometry of Modular Varieties}
\author[A. B. ~Goncharov ]{Alexander B. Goncharov}

\address{
Department of Mathematics,\br
Brown University,\br
Providence RI 02912,
USA} \email{sasha@math.brown.edu}

\begin{abstract}
We discuss two \emph{arithmetical} problems, at first glance
unrelated:

1) The properties of the multiple $\zeta$-values
\begin{equation} \label{lb*}
\zeta(n_{1},\dots ,n_{m}) : = \sum_{0 < k_{1} < k_{2} < \dots  < k_{m} }
\frac{1}{k_{1}^{n_{1}}k_{2}^{n_{2}}\dotsb k_{m}^{n_{m}}} \qquad n_m >1
\end{equation}
and their generalizations, multiple polylogarithms at $N$-th roots of unity.

2) The action of the absolute Galois group
on the pro-$l$ completion
$$
\pi^{(l)}_{1}(X_N):= \pi^{(l)}_{1}({\Bbb P}^{1} \backslash \{0, \mu_N,\infty  \}, v)
$$
 of the fundamental group of $X_N:= {\Bbb P}^{1} \backslash \{ 0, \infty$ and all $N$-th roots
of unity$\}$.

These problems are the Hodge and $l$-adic sites of the following one:

3) Study the Lie algebra of the image of
motivic Galois group acting
on the motivic fundamental group of
${\Bbb P}^{1} \backslash \{0, \mu_N,\infty\}$.

We will discuss   a
surprising connection between these problems
and \emph{geometry} of
the modular varieties
$$
Y_1(m:N):= \Gamma_1(m; N)\backslash GL_m(\R)/O_m\cdot \R^*
$$
 where
$\Gamma_1(m; N)$ is the subgroup of $GL_m(\Z)$ stabilizing $(0,\dots ,0,1)$ mod $N$.

In particular using this relationship we get precise results about
 the Lie algebra of the
image of the absolute Galois group in
${\rm Aut} \pi^{(l)}_{1}(X_N)$, and sharp estimates on the
dimensions of the $\Q$-vector spaces generated by the multiple polylogarithms
at $N$-th roots of unity, \emph{depth} $m$ and  \emph{weight} $w:= n_1+ \dots  +n_m$.

The simplest case of the problem 3) is related to  the classical theory of cyclotomic units.
Thus the subject of this lecture is \emph{higher cyclotomy}.
\end{abstract}

\maketitle

\section{The Multiple \boldmath$\zeta$-Values}

\subsection{The algebra of multiple \boldmath$\zeta$-values and its conjectural description}

Multiple $\zeta$-values (\ref{lb*})
were invented by L.~Euler \cite{E}. Euler discovered that the numbers
$\zeta(m,n)$, when $w:= m+n$  is odd, are $\Q$-linear combinations of
$\zeta(w)$ and $\zeta(k) \zeta(w-k)$. Then these numbers were neglected.
About $10$ years ago they were resurrected as the coefficients of Drinfeld's associator
\cite{Dr},
rediscovered by D.~Zagier \cite{Z}, appeared in works of
M.~Kontsevich  on knot invariants and   the author \cite{G1,G2}
on mixed Tate motives over
Spec~$(\Z)$. More  recently they showed up in quantum field theory \cite{Kr,Br},
deformation quantization and so on.

We say that $w := n_1+\dots +n_m$ is the \emph{weight}  and $m$
is the \emph{depth} of (\ref{lb*}).

Let  ${\mathcal Z}$   be the space
of $\Bbb Q$-linear combinations of multiple $\zeta$'s. It
is a  commutative algebra
 over $\Bbb Q$.  For instance
\begin{equation} \label{4-15.3}
  \zeta (m)\cdot  \zeta (n) =   \zeta (m,n) + \zeta(n,m) +  \zeta (m+n)
\end{equation}
because
\begin{equation} \label{4-15.4}
    \sum_{k_1,k_2 >0}\frac{1}{k_1^{ m}k_2^{n}} =
    \Bigl(\sum_{0 < k_1 < k_2 } + \sum_{0 < k_1 = k_2 }
    + \sum_{ k_1  > k_2 > 0}\Bigr)\frac{1}{k_1^{ m}k_2^{n}}\,.
\end{equation}

The  only known results about the classical $\zeta$-values are the following:
\begin{equation} \label{4-15.1}
\zeta(2n) = (-1)^{n-1}(2\pi)^{2n}\cdot\frac{B_{2n}}{2\cdot (2n)!}
\quad \mbox{(Euler)};\qquad \zeta(3) \not \in \Q  \quad \mbox{(Apery)}\,.
\end{equation}
Here $B_k$ are the Bernoulli numbers: $\frac{t}{e^t-1} = \sum B_kt^k/k!$.

To describe the hypothetical structure of
the algebra ${\mathcal Z}$ we introduce a free graded Lie algebra
${\mathcal F}(3,5,\dotsc )_{\bullet}$, which is freely
generated by elements $e_{2n+1}$ of degree $-(2n+1)$ where $n \geq 1$.
Let
$$
U{\mathcal F}(3,5,\dotsc)_{\bullet}^{\vee}:= \oplus_{n \geq 1}
\Bigl(U{\mathcal F}(3,5,\dotsc)_{-(2n+1)}\Bigr)^{\vee}
$$
be the graded dual to its universal
enveloping algebra. It is  $\Z_{+}$-graded.

\begin{conjecture} \label{1.25.1}
 a) The  weight provides a grading on the algebra ${\mathcal
Z}$.

b) One has an isomorphism of graded algebras over $\Q$
\begin{equation} \label{4-15.2}
{\mathcal Z}_{\bullet} = \Q[\pi^2] \otimes_{\Q} U{\mathcal F}(3,5,\dotsc)_{\bullet}^{\vee}
\qquad \quad\deg \pi^2 :=  2\,.
\end{equation}
\end{conjecture}

Part a) means that relations between $\zeta$'s of different weight, like
$\zeta(5) = \lambda \cdot \zeta(7)$ where $\lambda \in \Q$, are impossible.
For motivic interpretation/formulation of conjecture~\ref{1.25.1} see
Section~12 in \cite{G1}. For its $l$-adic version see conjecture~\ref{l-adic} below.

\begin{theorem} \label {th1.1}
One has $ \dim  {\mathcal Z}_{k}  \leq  \dim (\Q[\pi^2]
\otimes U{\mathcal F}(3,5,\dotsc)^{\vee})_{k}$.
\end{theorem}

Both the origin of conjecture~\ref{1.25.1} and proof of this theorem are
based on theory of mixed Tate motives
over $Spec(\Z)$: multiple $\zeta$-values are periods of framed mixed Tate motives
over $\Z$ (see \cite{G8}),
and one can \emph{prove} that the framed mixed Tate motives
over $\Z$ form an algebra which is isomorphic to the one appearing on the right hand side of (\ref{4-15.2}).
This gives theorem~\ref{th1.1}.
Conjecture~\ref{1.25.1} just means that \emph{every} such a period is given by multiple
$\zeta$-values.

For the definition of the {\it abelian}
category of mixed Tate motives over a number field
convenient for our approach see chapter~5 in \cite{G6}. It has
all the expected properties and based on V.~Voevodsky's
construction of the triangulated category of motives
\cite{V}. Another approach to mixed motives has been developed by
M.~Levine \cite{L1}, \cite{L}. A construction of the framed mixed
Tate motive over $\Q$ related to multiple $\zeta$'s can be obtained
by combining constructions in section~12 of \cite{G2} and  chapter~5 of
\cite{G6}.

It is difficult to estimate ${\rm dim}  {\mathcal Z}_{k}$ from below: we
believe that ${ \zeta}(5) \not \in \Q$ but nobody can prove it.

One may reformulate  conjecture~\ref{1.25.1} as a hypothetical description of the $\Q$\guio{vector}
space ${\mathcal P}{{\mathcal Z}}$ of primitive multiple $\zeta$'s:
\begin{equation} \label{lbb}
{\mathcal P}{{\mathcal Z}}_{ \bullet} :=  \frac{  {\mathcal Z}_{\bullet}}
 { {\mathcal Z }_{>0} \cdot
 {\mathcal Z}_{>0}} \stackrel{?}{=} <\pi^2> \oplus {\mathcal F}(3,5,\dotsc)_{\bullet}^{\vee}\,.
\end{equation}
Here $<\pi^2>$ is a $1$-dimensional $\Q$-vector space generated by $\pi^2$, and
${\mathcal Z }_{>0}$  is generated by $\pi^2$
and $\zeta(n_1, \dots , n_m)$.

\begin{example}
There are $2^{10}$ convergent multiple $\zeta$'s
of the weight $12$. However according to theorem~\ref{th1.1} ${\rm dim} {\mathcal
Z}_{12} \leq 12$. One should have   ${\rm dim} {\mathcal P}{\mathcal  Z}_{ 12} = 2$
     since ${\mathcal F}(3,5,\dotsc)_{-12}$ is spanned
over  $\Q$
by  $[e_{5},e_{7}]$ and $[e_{3},e_{9}]$.
The $\Q$-vector space of  \emph{decomposable} multiple $\zeta$'s of the weight $12$
is supposed to be generated by
$$
\pi^6, \quad  \pi^3{\zeta}(3)^2,  \quad
\pi^2 {\zeta}(3)\zeta(5), \quad \pi^2 {\zeta}(3,5), \quad \pi {\zeta}(3)\zeta(7),
\quad \pi {\zeta}(5)^2, \quad \pi {\zeta}(3,7),
$$
$$
\quad {\zeta}(3)^4, \quad {\zeta}(5){\zeta}(7),
\quad
{  \zeta}(3)\zeta(9)\,.
$$
\end{example}

The algebra  $U{\mathcal F}(3,5,\dotsc)_{\bullet}^{\vee}$ is commutative.
It is isomorphic to the space of
noncommutative polynomials in variables $f_{2n+1}$, $n=1,2,3,\dotsc$
with the algebra structure given by the shuffle product.

Let ${\mathcal F}(2,3)_{\bullet}$ be the free graded Lie algebra generated by two
elements of degree $-2$ and $-3$.
Its graded dual $U{\mathcal F}(2,3)_{\bullet}^{\vee}$ is isomorphic as a
graded vector space  to the space of noncommutative
polynomials in two variables $p$ and $g_3$ of degrees 2 and 3.
There is canonical isomorphism of \emph{graded vector spaces}
$$
\Q[\pi^2] \otimes U{\mathcal F}(3,5,\dotsc
)_{\bullet}^{\vee} = U{\mathcal F}(2,3)_{\bullet}^{\vee}\,.
$$
The rule is clear from the pattern
$
(\pi^2)^3 f_3(f_7)^3(f_5)^2 \longrightarrow p^3 g_3(g_3p^2)^3(g_3p)^2
$.

In particular if $d_k := {\rm dim} {\mathcal Z}_{k}$ then one should have $d_k =
d_{k-2} + d_{k-3}$. This rule has been observed in computer
calculations of D.~Zagier
for $k\leq 12$. Later on
extensive computer calculations, confirming it, were made by D.~Broadhurst~\cite{Br}.

\subsection{The depth filtration}

Conjecture~\ref{1.25.1}, if true, would give a very simple and clear picture for
the structure of the multiple $\zeta$-values algebra.
However this algebra has an additional structure: the depth filtration, and conjecture
\ref{1.25.1} tells us nothing about it.
The study of the depth filtration moved the subject in a completely unexpected direction: towards
geometry of modular varieties for $GL_m$.

To formulate some results about the depth filtration consider the algebra $\overline {\mathcal Z}$
spanned over $\Q$ by the numbers
$$
\overline \zeta(n_1, \dots , n_m):= (2\pi i)^{-w}\zeta(n_1, \dots , n_m)\,.
$$
It is filtered by the weight and depth. Since
$\overline \zeta(2) = -1/24$, there is no weight grading anymore.
Let ${\rm Gr}^{W,D}_{w,m} {\mathcal P}\overline {\mathcal Z}$ be the associated graded.
We assume that $1$ is of depth $0$.
Denote by $d_{w,m}$ its dimension over $\Q$.

Euler's  classical computation of $\zeta(2n)$ (see~(\ref{4-15.1})) tells us that
$d_{2n,1} =0$. Generalizing this it is not hard to prove that $d_{w,m}=0$
if $w+m$ is odd.

\begin{theorem} \label{1.25.2}
\begin{enumerate}[a)]
\item
\begin{equation} \label{1.25.2+}
 d_{w,2}  \leq    \left[\frac{w-2}{6}\right] \quad
\mbox{if $w$ is even}\,.
\end{equation}
\item
$$
  d_{w,3}  \leq  \left[\frac{(w-3)^2-1}{48}
    \right] \quad \mbox{if $w$ is odd}\,.
$$
\end{enumerate}
\end{theorem}
The part a) is due to Zagier;
the dimension of the space of cusp forms for $SL_2(\Bbb Z)$ showed up in
his investigation
of the double shuffle relations for the depth two multiple $\zeta$'s, (\cite{Z}).
The part b) has been proved in \cite{G4}. Moreover
we proved that, assuming  some standard conjectures in
arithmetic algebraic geometry,  these estimates are exact,
see also corollary~\ref{4-30.1} and theorem~\ref{mth1t}.

\begin{problem} \label{4-16.2} Define explicitly a depth filtration on
the Lie 
coalgebra
${\mathcal F}
(3,5,\dotsc)^{\vee}$ 
which under the isomorphism (\ref{lbb})
should correspond to the
depth filtration on the space of primitive multiple $\zeta$-values.
\end{problem}

The cogenerators of the Lie coalgebra ${\mathcal F}(3,5,\dotsc )^{\vee}$ correspond to
$\zeta(2n+1)$. So  a naive guess would be that the dual to the
lower central series filtration on
${\mathcal F}(3,5,\dotsc)$ coincides with the  depth filtration.
However then one should have
$d_{12,2} =2$, while according to formula (\ref{1.25.2+}) $d_{12,2} =1$.
Nevertheless ${\rm dim}{\mathcal P}{\mathcal Z}_{12} =2$, but the new
transcendental number appears only in the  depth $4$.

\subsection{A heuristic discussion}

Conjecture~\ref{1.25.1} in the form (\ref{lbb}) tells us that
the space of primitive
multiple $\zeta$'s should have a Lie coalgebra structure.
How to determine its coproduct $\delta$
in terms of the multiple $\zeta$'s?
Here is the answer for the depth $1$ and $2$ cases.
(The general case later on).
Consider the generating series
$$
\zeta(t):= \sum_{m>0} \zeta(m)t^{m-1}, \qquad  \zeta(t_1,t_2):  =
\sum_{m,n >0} \zeta(m,n)t_1^{m-1}t_2^{n-1}
$$
Then $\delta \zeta(t) =0$, i.e.\ $\delta \zeta(n) =0$ for all $n$, and
\begin{equation} \label{4-30.4}
\delta \zeta(t_1,t_2) = \zeta(t_2) \wedge \zeta(t_1) + \zeta(t_1) \wedge
\zeta(t_2-t_1)  - \zeta(t_2) \wedge
\zeta(t_1-t_2)\,.
\end{equation}
To make sense out of this
 we have to go from  the numbers
$\zeta(n_1, \dots , n_m)$ to their more structured
counterparts:
framed mixed Tate motives $\zeta_{{\mathcal M}}(n_1, \dots , n_m)$,
or their Hodge or $l$-adic realizations, (see \cite{G3} or Section~12 of \cite{G1}).
The advantage is immediately seen: the coproduct $\delta_{{\mathcal M}}$
is well defined by the general formalism (see~Section~10 in \cite{G1}),
one easily proves not only that $\zeta_{{\mathcal M}}(2n) =0$
(motivic version of Euler's theorem) as well as
$\zeta_{{\mathcal M}}(1) =0$, but also that
$\zeta_{{\mathcal M}}(2n+1) \not =0$, and there are no linear relations between
$\zeta_{{\mathcal M}}(2n+1)$'s!  Hypothetically we loose no information:
linear relations between the multiple $\zeta$'s should reflect  linear relations
between their motivic avatars.
Using $\zeta_{{\mathcal M}}(2n) =0$
we rewrite formula (\ref{4-30.4}) as
\begin{equation} \label{4-30.30}
\delta_{{\mathcal M}}: \zeta_{{\mathcal M}}(t_1,t_2) \lms  \Bigl(1 + U + U^2 \Bigr)
\zeta_{{\mathcal M}}(t_1) \wedge \zeta_{{\mathcal M}}(t_2)
\end{equation}
where $U$ is the linear operator $(t_1, t_2) \lms (t_1-t_2, t_1)$.
For example $\delta_{{\mathcal M}}$ sends  the subspace of
weight $12$ double $\zeta_{{\mathcal M}}$'s to
 a one dimensional $\Q$-vector space generated by
$3 \zeta_{{\mathcal M}}(3) \wedge \zeta_{{\mathcal M}}(9) +
\zeta_{{\mathcal M}}(5) \wedge \zeta_{{\mathcal M}}(7)$.
One can identify the cokernel of the map (\ref{4-30.30}),
restricted to the weight $w$ subspace,
  with $H^1(GL_2(\Z), S^{w-2}V_2 \otimes \varepsilon_2)$ where
$V_2$ is
the standard $GL_2$-module, and $ \otimes \varepsilon_2$
is the twist by the determinant,
i.e. with the
space of weight $w$ cusp forms for $GL_2(\Z)$.
Moreover, one can prove that  ${\rm Ker}\delta_{{\mathcal M}}$ is spanned by
$\zeta_{{\mathcal M}}(2n+1)$'s: this is a much more difficult result which uses all the machinery of mixed motives. Thus an element of the depth $2$ associated graded
of the space of primitive double $\zeta$'s is zero
if and only if its coproduct is $0$. So
formula (\ref{4-30.30}) provides a complete
description of the space of double $\zeta$'s. In particular $d_{12,2} =1$.

For the rest of this paper we suppress the motives working
 mostly with the $l$-adic side of the  story and  looking
at the Hodge side for motivations.

\section{Galois Symmetries of the pro-\boldmath$l$ 
Completion of the Fundamental
Group of \boldmath$ {\Bbb P}^1 \backslash \{0, \mu_N, \infty\}$}

\subsection{The Lie algebra of the image of the Galois group}

Let $X$ be a regular curve, $\overline X$ the corresponding projective curve,
and $v$ a tangent vector at a point $x \in \overline X$.
According to Deligne \cite{D} one can define the geometric
profinite fundamental group $\widehat \pi_1(X, v)$ based at $v$.
If $X$, $x$ and $v$ are defined over a  number field $F$
then the group ${\rm Gal}_F:= {\rm Gal}(\overline  \Q/F) $
acts by automorphisms of $\widehat \pi_1(X, v)$.

If $X = {\Bbb P}^{1} \backslash \{0, \mu_N,\infty \}$ there is a tangent vector $v_\infty$
corresponding to the inverse $t^{-1}$ of the canonical coordinate $t$ on
${\Bbb P}^{1} \backslash \{0, \mu_N,\infty \}$.
Denote by $\pi^{(l)}$ the pro-$l$-completion of the group
$\pi$.
We will investigate the map
\begin{equation} \label{121212**}
\Phi^{(l)}_{N}: {\rm Gal}_{\Bbb Q} \longrightarrow
{\rm Aut} \pi^{(l)}_{1}({\Bbb P}^{1} \backslash \{0, \mu_N,\infty \}, v_\infty )\,.
\end{equation}
When $N=1$ it was studied by Grothendieck \cite{Gr}, Deligne \cite{D},
Ihara (see \cite{Ih1,Ih3}), Drinfeld \cite{Dr}, and others
(see~\cite{Ih2}), but for $N>1$ it was not investigated.

Denote by $H(m)$ the lower central series for the group
$H$. Then the quotient
$\pi^{(l)}_{1}(X_N)/ \pi^{(l)}_{1}(X_N)(m)
$ is an $l$-adic Lie group.
Taking its Lie algebra and making the projective limit
over $m$ we get a pronilpotent Lie algebra over $\Q_l$:
$$
{\Bbb L}^{(l)}_N:= \lim_{\longleftarrow}Lie \Bigl( \frac{\pi^{(l)}_{1}(X_N)}
{ \pi^{(l)}_{1}(X_N)(m)}\Bigr)\,.
$$
Similarly one defines an $l$-adic pronilpotent Lie algebra ${\Bbb L}^{(l)}(X, v)$
corresponding to the geometric fundamental group $\widehat \pi_1(X,v)$ of a
variety $X$ with a base at $v$.

For the  topological reasons ${\Bbb L}^{(l)}_N$
is a free pronilpotent Lie algebra over $\Q_l$
with $n+1$ generators corresponding to the loops around $0$ and $N$-th roots of unity.

Let $\Q(\zeta_n)$ be the field generated by $n$-th roots of unity. Set
$\Q(\zeta_{l^{\infty}N}):= \cup \Q(\zeta_{l^{a}N})$.
We restrict map (\ref{121212**}) to the Galois group
${\rm Gal}_{\Q(\zeta_{l^{\infty}N})}$. Passing to Lie algebras we get a homomorphism
$$
\phi^{(l)}_{N}\colon {\rm Gal}_{\Q(\zeta_{l^{\infty}N})} \lra {\rm Aut}({\Bbb L}^{(l)}_N)\,.
$$
Let us linearize  the image of this map.
Let ${\Bbb L}^{(l)}_N(m)$ be the lower central series for
the Lie algebra ${\Bbb L}^{(l)}_N$. There are homomorphisms to $l$-adic Lie groups
$$
\phi^{(l)}_{N;m}\colon {\rm Gal}_{\Q(\zeta_{l^{\infty}})} \lra {\rm Aut}
\Bigl({\rm L}^{(l)}_N/{\Bbb L}^{(l)}_N(m)\Bigr)\,.
$$
The main hero of this story is the pronilpotent Lie algebra
$$
{\mathcal G}^{(l)}_{N}:= \lim_{\stackrel{\longleftarrow}{m}} Lie
\Bigl(Im \phi^{(l)}_{N;m}\Bigr) \hookrightarrow {\rm Der} {\Bbb L}^{(l)}_N\,.
$$
When $N=1$ we denote it by ${\mathcal G}^{(l)}$.

\begin{conjecture} \label{l-adic}
${\mathcal G}^{(l)}$ is a free
Lie algebra with generators indexed by odd integers $\geq 3$.
\end{conjecture}

It has been formulated, as a question, by Deligne \cite{D} and Drinfeld \cite{Dr}.

\subsection{The weight and depth filtration on \boldmath${\Bbb L}^{(l)}_N$}

There are two increasing
filtrations by ideals on the
Lie algebra ${\Bbb L}^{(l)}_N$, indexed by negative integers.

\emph{The  weight filtration ${\mathcal F}^{W}_{\bullet}$}.
It coincides with the lower central series for  ${\Bbb L}^{(l)}_N$:
$$
{\Bbb L}^{(l)}_N = {\mathcal F}^{W}_{  -1}{\Bbb L}^{(l)}_N; \qquad {\mathcal F}^{W}_{ -n -1}{\Bbb L}^{(l)}_N :=
[{\mathcal F}^{W}_{-n}{\Bbb L}^{(l)}_N, {\Bbb L}^{(l)}_N]\,.
$$

\emph{The depth filtration ${\mathcal F}^{D}_{\bullet}$}. The natural inclusion
$$
{\Bbb P}^{1} \backslash \{0, \mu_N,\infty \} \hookrightarrow
{\Bbb P}^{1} \backslash \{0, \infty \}
$$
provides a morphism of the corresponding fundamental Lie algebras
$$
p\colon {\Bbb L}^{(l)}_N \lra {\Bbb L}^{(l)}({\Bbb P}^{1} \backslash \{0, \infty \}) = \Q_l(1)\,.
$$
Let   ${\mathcal I}_N$
be the kernel of this projection.
Its powers   give
 the  depth filtration:
$$
{\mathcal F}^{D}_{ 0}{\Bbb L}^{(l)}_N = {\Bbb L}^{(l)}_N, \quad
{\mathcal F}^{D}_{ -1}{\Bbb L}^{(l)}_N = {\mathcal I}_N,
\quad {\mathcal F}^{D}_{ -n-1}{\Bbb L}^{(l)}_N = [{\mathcal I}_N,
{\mathcal F}^{D}_{ -n}{\Bbb L}^{(l)}_N]\,.
$$

\subsection{The Galois Lie algebra and its shape}
 
These filtrations induce  two filtrations on the Lie algebra
${\rm Der}{\Bbb L}^{(l)}_N$ and hence on the Lie algebra ${\mathcal G}^{(l)}_{N} $.
The  associated graded Lie algebra  ${\rm Gr}{\mathcal G}_{\bullet \bullet}^{(l)}(\mu_N)$, which we call the level $N$ Galois Lie algebra,  is
 bigraded
by the weight $-w$ and depth $- m $.
The weight filtration can be  defined by a grading. Moreover one can define it in a way
compatible the depth filtration and the subspace ${\mathcal G}_N^{(l)}
\subset {\rm Der}{\Bbb L}^{(l)}_N$. Therefore
$$
{\rm Gr}{\mathcal G}_{\bullet \bullet}^{(l)}(\mu_N) \hookrightarrow
{\rm Gr}_{\bullet \bullet}{\rm Der}{\Bbb L}_N^{(l)}\,.
$$

{\it The depth $m$ quotients}. For any $m\geq 1$ there is the depth $\geq -m$ 
(we will also  say depth $m$) quotient Galois Lie algebra:
\begin{equation} \label{grkv}
{\rm Gr}{\mathcal G}_{\bullet, \geq -m}^{(l)}(\mu_N) :=
\frac{{\rm Gr}{\mathcal G}_{\bullet, \bullet}^{(l)}(\mu_N) }
{{\rm Gr}{\mathcal G}_{\bullet, < -m}^{(l)} (\mu_N)}\,.
\end{equation}
It is a nilpotent graded Lie algebra of the nilpotence class $\leq m$. 

{\it The diagonal Lie algebra}. Notice that ${\rm Gr}{\mathcal G}_{-w, -m}^{(l)}(\mu_N) =0$ if $w < m$. 
We define the diagonal Galois Lie algebra
${\rm Gr}{\mathcal G}^{(l)}_{\bullet}(\mu_N)$
as the Lie subalgebra of
${\rm Gr}{\mathcal G}^{(l)}_{\bullet \bullet}(\mu_N)$ formed by the 
 components
with $w = m$. It is graded by the weight. It 
can be defined as the Lie subalgebra of 
${\rm Gr}^W{\mathcal G}^{(l)}(\mu_N)$ by 
imposing the $weight = depth$ condition. Thus we do not need 
to take the associated graded  for the depth filtration for its 
definition. Therefore the diagonal Galois Lie algebra is isomorphic, although  
non canonically, to a Lie subalgebra of ${\mathcal G}^{(l)}(\mu_N)$.

The picture below  exibites all possibly  non zero components of the 
Galois Lie algebra and indicates its depth 2 quotient and the diagonal Lie subalgebra. 

\begin{center}
\hspace{4.0cm}
\epsffile{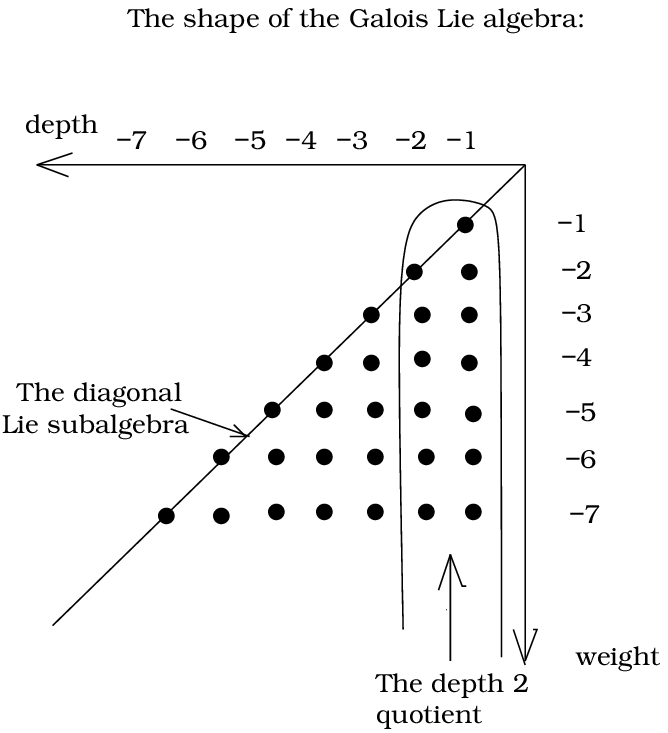}
\end{center}

\subsection{The mysterious correspondence}

Let $V_m$ be the standard $m$-dimensional representation of $GL_m$.
Our key point (\cite{G2}-\cite{G5}) is that

\vspace{3mm}
\noindent
\begin{minipage}[c]{0.3\textwidth}
the structure of Galois Lie algebra
$$
{\rm Gr}{\mathcal G}_{\bullet, \geq -m}^{(l)}(\mu_N)$$
\end{minipage}
\hspace{4mm}\begin{minipage}[c]{0.15\textwidth}
\emph{is related to}
\end{minipage}\hspace{4mm}
\begin{minipage}[c]{0.4\textwidth}
geometry of local systems with fibers $S^{\bullet-m}{\rm V}_m$ over (the closure of)
modular variety $Y_1(m;N)$, which is defined for $m>1$ as
$$
\Gamma_1(m;N)\backslash GL_m(\R)/O_m \cdot {\R}^*\,.
$$
\end{minipage}
\vspace{5mm}

\noindent
The adelic approach to modular varieties shows that for $m=1$ we have 
\begin{equation} \label{7.24.00.1}
 Y_1(1; N):= \quad S_N:= \quad {\rm Spec}\Z[\zeta_N][\frac{1}{N}]
\end{equation}

In particular the diagonal level $N$ Galois Lie algebra is 
related to the geometry of the modular varieties $Y_1(m; N)$.

Both the dual to the Galois Lie algebras (\ref{grkv}) and
the modular varieties ${\overline Y}_1(m;N)$ form inductive systems with respect to $m$.
The correspondence  is compatible with these inductive structures.

Recall the standard cochain complex of a Lie algebra ${\mathcal G}$
$$
{\mathcal G}^{\vee} \stackrel{\delta}{\lra} \Lambda^2 {\mathcal G}^{\vee}
\stackrel{\delta}{\lra} \Lambda^3 {\mathcal G}^{\vee} \lra \dotsb
$$
where the first differential is dual to the commutator map $[ , ]\colon \Lambda^2 {\mathcal G} \lra {\mathcal G}$,
and the others are obtained using
the Leibniz rule. The condition $\delta^2 =0$ is equivalent to the Jacobi identity.

For a precise form of this correspondence see section~6. It relates the depth
$m$, weight $w$ part of the standard cochain complex of the Lie algebra (\ref{grkv})
with
\begin{equation} \label{4-26.1}
\mbox{(rank $m$ modular complex) $\otimes_{\Gamma_1(m;N)}S^{w-m}V_m$}\,.
\end{equation}
The rank $m$ modular complex is a complex of $GL_m(\Z)$-modules constructed
purely combinatorially.
It has a \emph{geometric realization} in the symmetric space
${\Bbb H}_m := GL_m(\R)/O(m) \cdot \R^*$,
see section~7 and \cite{G9}. It is well understood only for $m \leq 4$.

\subsection{Examples of this correspondence}

Strangely enough it is more convenient to describe the structure of 
the bigraded Lie algebra
$$
{\rm Gr}\widehat {\mathcal G}_{\bullet \bullet}^{(l)}(\mu_N):=
{\rm Gr}{\mathcal G}_{\bullet \bullet}^{(l)}(\mu_N)\oplus {\Q_l}(-1,-1)
$$
where ${\Q_l}(-1,-1)$ is a one dimensional Lie algebra of weight and   depth $-1$.
Its motivic or Galois-theoretic meaning is non clear: it should correspond
to  $\zeta(1)$.

\textbf{a) The depth \boldmath$1$ case}. The depth $-1$ quotient
${\rm Gr}_{\bullet, -1}{\mathcal G}^{(l)}_N$is an abelian Lie algebra.
Its structure is described by the following theorem.

\begin{theorem} \label{thsdi}There is a natural isomorphism of $\Q_l$-vector spaces
\begin{equation} \label{demfed+}
{\rm Hom}\Bigl(K_{2n-1}(\Z[\zeta_N, N^{-1}]), \Q_l\Bigr) \quad \stackrel{=}{\lra} \quad
{\rm Gr}_{-n, -1}{\mathcal G}_N^{(l)}\,.
\end{equation}
\end{theorem}

According to the Borel theorem one has
\begin{equation} \label{demfed}
{\rm dim} K_{2n-1}(\Z) \otimes \Q= \left\{ \begin{array}{ll}
0 &  \quad \mbox{$n$: even} \\
 1 &    \quad \mbox{$n>1$: odd} \end{array} \right.
 \end{equation}
\begin{equation} \label{demfed&&}
{\rm dim} K_{2n-1}(\Z[\zeta_N, N^{-1}]) \otimes \Q =  \left\{ \begin{array}{ll}
\frac{\varphi(N)}{2} &  \quad \mbox{$N>2, n>1$} \\
\frac{\varphi(N)}{2}+ p(N)-1 &  \quad \mbox{$N>2, n=1$} \\
 1 &    \quad \mbox{$N=2$} \end{array} \right.\,.
 \end{equation}
where $p(N)$ is the number of prime factors of $N$.

Theorem~\ref{thsdi} for $N=1$  is known thanks to Soul\'e, Deligne \cite{D}, and
Ihara \cite{Ih2}. The general case can be deduced from the motivic theory of classical
polylogarithms developed by Deligne and  Beilinson \cite{D,BD}.
In the case $n=1$  there is canonical isomorphism justifying the name
``higher cyclotomy'' for our story:
\begin{equation}\label{4-30.2}
{\rm Gr}{\mathcal G}_{-1, -1}^{(l)}(\mu_N) =
{\rm Hom}\Bigl(\mbox{group of the cyclotomic units in $
\Z[\zeta_N][\frac{1}{N}]$},  \quad \Q_l\Bigr)\,.
\end{equation}

The level $N$ modular variety for ${GL_1}_{/{\Q}}$ is the  scheme
$S_N$, see (\ref{7.24.00.1}). It has $\varphi(N)$ complex points
parametrized by the primitive roots of unity
$\zeta_N^{\alpha}$ where $(\alpha, N)=1$. Our correspondence for  $m=1$ is given by
the isomorphism (where $+$ means invariants under the complex conjugation):
$$
[\Q(n-1)-\mbox{valued functions on} \quad S_N\otimes\C]^+ \quad
\stackrel{=}{\lra} \quad K_{2n-1}(S_N)\otimes \Q
$$
provided by motivic classical polylogarithms: one associates to $\zeta_N^{\alpha}$ the the cyclotomic element
$\{\zeta_N^{\alpha}\}_n \in K_{2n-1}(S_N)$, whose regulator is computed via
$Li_n(\zeta_N^{\alpha})$.

\textbf{b) The depth \boldmath$2$ case, \boldmath$N=1$}. The structure of the depth
$\geq -2$ quotient of the Lie algebra
${\rm Gr}\widehat {\mathcal G}_{\bullet, \bullet}^{(l)} $
is completely described by the commutator map
\begin{equation} \label{dep2}
[,]\colon \quad \Lambda^2 {\rm Gr}\widehat {\mathcal G}_{-1, \bullet}^{(l)} \lra {\rm Gr}
{\mathcal G}_{-2, \bullet}^{(l)}\,.
\end{equation}

\emph{Construction of the dual to complex (\ref{dep2})}. Look at
the classical modular triangulation of the hyperbolic plane ${\Bbb H}_2$ where
the central ideal triangle has vertices at $0, 1, \infty$:
\begin{center}
\epsffile{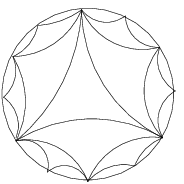}
\end{center}
The group $GL_2(\R)$,  acting on $\C \backslash \R$ by $z \lms
\frac{az+b}{cz+d}$
commutes, with $z \lms \overline z$. We let $GL_2(\R)$ act on ${\Bbb
H}_2$ by identifying
${\Bbb H}_2$ with the quotient of $\C \backslash \R$ by complex
conjugation. The subgroup $GL_2(\Z)$
preserves the modular picture.
Consider the chain complex of the modular triangulation placed in
degrees $[1,2]$:
\begin{equation} \label{dep44}
M^*_{(2)}:= \quad M^1_{(2)} \lra M^2_{(2)}\,.
\end{equation}
It is a complex of $GL_2(\Z)$-modules.
The group $M^1_{(2)}$ is generated by the triangles, and $M^2_{(2)}$ by the geodesics.
Let $\varepsilon_2$ be the one dimensional $GL_2$-module given by the determinant.
\begin{lemma} \label{modcoho}
Let $\Gamma$ be a finite index subgroup of $GL_2(\Z)$ and $V$ a
$GL_2$-module over $\Q$. Then the complex
 ${M}_{(2)}^{\ast} \otimes_{\Gamma} V [1]$
computes the cohomology $H^{*}(\Gamma,  V\otimes \varepsilon_2)$.
\end{lemma}
\begin{theorem} \label{dw1}
The weight $w$ part of the dual to complex (\ref{dep2})
is \textbf{canonically} isomorphic to the complex
\begin{equation} \label{dep4}
\Bigl(M^*_{(2)} \otimes_{GL_2(\Z)} S^{w-2}{\rm V}_2\Bigr)\otimes \Q_l\,.
\end{equation}
\end{theorem}

Motivic version of this theorem was obtained in section~7 of \cite{G2}.
Its Hodge side provides a refined version of the story told in section~1.3.

According to the lemma  complex (\ref{dep4}) computes
$
H^{*-1}(GL_2(\Z), S^{w-2}{\rm V}_2\otimes \varepsilon_2)
$.
Since these cohomology groups
  are known, we compute the Euler characteristic of the complex (\ref{dep2}) and
using theorem~\ref{thsdi} and  formula (\ref{demfed})`
get the following result, the $l$-adic version of
theorem~\ref{1.25.2}, (see related results of Ihara and Takao in \cite{Ih3}).
\begin{corollary} \label{4-30.1}
\begin{equation} \label{demfedd}
\qquad {\rm dim} {\rm Gr}{\mathcal G}_{-w, -2}^{(l)}  =
\quad    \left\{ \begin{array}{ll}
0 &  w: \quad \mbox{odd}  \\
   \left[ \frac{w-2}{6} \right]  &    w:  \quad \mbox{even}\,. \end{array} \right.
 \end{equation}
\end{corollary}

\textbf{c) \boldmath$N= p$ is a prime, \boldmath$w=m=-2$.}
The structure of the weight  $\geq -2$ 
quotient of the \emph{diagonal} Galois Lie algebra
${\rm Gr}{\mathcal G}^{(l)}_{\bullet}(\mu_N)$ 
is described by the commutator map
\begin{equation} \label{2/7/00/2}
[,]\colon  \Lambda^2{\rm Gr}\widehat {\mathcal G}_{-1, -1}^{(l)}(\mu_p) \lra
{\rm Gr}{\mathcal G}_{-2, -2}^{(l)}(\mu_p)
 \end{equation}
 Projecting the modular triangulation of the hyperbolic plane onto the modular curve
$Y_1(p):= \Gamma_1(p) \backslash {\Bbb H}_2$ we get  the modular triangulation of
$Y_1(p)$. The complex involution acts on the modular curve preserving the
triangulation. Consider the  following complex, where $+$ means invariants of the
complex involution:
\begin{equation} \label{MAN}
\Bigr(\mbox{{\rm the
chain complex of the modular triangulation of $Y_1(p)$} }\Bigl)^+ \otimes \Q_l\,.
\end{equation}

\begin{theorem} The dual to complex (\ref{2/7/00/2})
is naturally \textbf{isomorphic} to complex (\ref{MAN}).
\end{theorem}

In particular there is canonical isomorphism
\begin{equation} \label{MAN*}
\Q_l[\mbox{triangles of the modular triangulation of $Y_1(p)$}]^+ =
\Bigl( {\rm Gr}{\mathcal G}_{-2, -2}^{(l)}(\mu_p)\Bigr)^{\vee}\,.
\end{equation}
Computing the Euler characteristic of the complex (\ref{2/7/00/2}) using
(\ref{MAN*}) and (\ref{4-30.2}) we get
$$
{\rm dim} {\rm Gr}{\mathcal G}^{(l)}_{-2}(\mu_p) = \frac{(p-5)(p-1)}{12}\,.
$$

Deligne proved \cite{D2} that the  Hodge-theoretic version of ${\mathcal G}^{(l)}_{N}$
is free when $N=2$, and very recently extended the arguments to the case 
$N=3,4$. The results above imply that it can not be free for sufficiently big $N$.
For instance it is not free for a prime  $N=p$ if  the genus of $Y_1(p)$ is
positive, i.e.\ $p>5$. Indeed,
$$
 \mbox{the depth $\leq 2$ part of\ }  H^2({\mathcal G}^{(l)}_p) =
 H^2_{(2)}({\mathcal G}^{(l)}_{\bullet}(\mu_p)) =  H^1(\Gamma_1(p), \varepsilon_2)\,.
$$

\subsection{Our strategy}

To describe the structure of the Galois Lie algebras (\ref{grkv})
in general we need the  dihedral Lie algebra of the group $\mu_N$
(\cite{G4,G5}) recalled in section~4.
To motivate to some extent its definition we turn in section~3
 to the Hodge side of higher cyclotomy.
As explained in section~5 both Galois and dihedral Lie algebra of $\mu_N$
act in a special way on the pronilpotent completion of
$\pi_1({\Bbb P^1} \backslash \{0, \mu_N, \infty\}, v_{\infty})$, and the Galois is contained in
the dihedral Lie algebra (\cite{G5}). In section~6 we relate the standard cochain complex of the
dihedral Lie algebra of $\mu_N$ with the modular complex (\cite{G4}), whose
 canonical  geometric realization in the symmetric space (\cite{G9}) is given in section~7.
Thus we
 related the structure of the Galois Lie algebras with geometry of modular varieties.

\section{Multiple Polylogarithms and higher Cyclotomy}

\subsection{Definition and iterated integral presentation}

Multiple polylogarithms (\cite{G1,G2})
 are defined as the power series
\begin{equation} \label{zhe5}
Li_{n_{1},\dots ,n_{m}}(x_{1},\dots ,x_{m})
=  \sum_{0 < k_{1} < k_{2} < \dots  < k_{m} } \frac{x_{1}^{k_{1}}x_{2}^{k_{2}}
\dotsb  x_{m}^{k_{m}}}{k_{1}^{n_{1}}k_{2}^{n_{2}}\dotsb k_{m}^{n_{m}}}
\end{equation}
generalizing both the
classical polylogarithms $Li_n(x)$ (if $m=1$)
and  multiple $\zeta$-values (if $x_1 = \dots  = x_m =1$).
These series are convergent for $|x_i| <1$.

Recall a definition of iterated integrals.
Let $\omega_{1},\dots ,\omega_{n}$ be $1$-forms on a manifold $M$ and
$\gamma \colon [0,1] \rightarrow M$  a path.
The iterated integral $\int_{\gamma} \omega_{1} \circ \dots  \circ \omega_{n}$
is defined inductively:
\begin{equation} \label{h1}
\int_{\gamma} \omega_{1} \circ \dots  \circ \omega_{n} \quad : = \quad \int_{0}^{1}(\int
_{\gamma_{t}} \omega_{1} \circ \dots  \circ \omega_{n-1}) \gamma^{*} \omega_{n}\,.
\end{equation}
Here $ \gamma_{t}$ is the restriction of $\gamma$ to the interval $[0,t]$
and
$\int_{\gamma_{t}} \omega_{1} \circ \dots  \circ \omega_
{n-1}$ is considered as a  function on $[0,1]$. We multiply it by the
 $1$-form $\gamma_t^{*} \omega_{n}$  and integrate.

Denote by  $I_{n_1,\dots ,n_m}(a_1:\dots :a_m:a_{m+1})$ the iterated integral
\begin{equation} \label{2/7/00/1}
\int_{0}^{a_{m+1}} \underbrace
 {\frac{dt}{ a_{1} - t} \circ \frac{dt}{t} \circ \dots  \circ
\frac{dt}{t}}_{n_{1} \quad  \mbox {times}} \circ   \dots   \circ
\underbrace {\frac{dt}{ a_{m} - t} \circ \frac{dt}{t} \circ \dots  \circ \frac{dt}{t}\,.}_
{n_{m} \quad  \mbox {times}}
\end{equation}
Its value depends only on the homotopy class of a path connecting $0$ and $a_{m+1}$ on $\C^* \backslash \{a_1, \dots , a_{m}\}$.
Thus it is a multivalued analytic function of $a_1, \dots , a_{m+1}$.
The following result provides an analytic continuation of multiple polylogarithms.

\begin{theorem} \label{BB}
 $Li_{n_1,\dots ,n_m}(x_1,\dots ,x_m) = I_{n_1,\dots ,n_m}(1: x_1: x_1x_2: \dots  : x_1\dotsb x_m)$.
\end{theorem}
The proof is easy:
develop $dt/(a_i - t)$ into a geometric series and integrate.
If $x_i =1$ we get the Kontsevich formula. In particular in the depth one case we recover the classical Leibniz presentation for $\zeta(n)$:
\begin{equation} \label{5*5}
{  \zeta}(n)  =
 \int_{0}^{1} \underbrace {\frac{dt}{1-t}  \circ \frac{dt}{t} \circ
\dots  \circ
\frac{dt}{t}\,. }
_{ n \quad  \mbox {times}}
\end{equation}

\subsection{Multiple polylogarithms at roots of unity}

Let $\overline {\mathcal Z}_{\leq w}(N)$ be the $\Q$-vector space spanned by the numbers
\begin{equation} \label{mLv}
\overline Li_{n_1,\dots ,n_m}( \zeta_N^{\alpha_1},\dots ,\zeta_N^{\alpha_m}):=
(2 \pi i)^{-w} Li_{n_1,\dots ,n_m}( \zeta_N^{\alpha_1},\dots ,\zeta_N^{\alpha_m});
\qquad  \zeta_N:= e^{2\pi i/N}\,.
\end{equation}
Here we may take any branch of $Li_{n_1,\dots ,n_m}(x_1,\dots ,x_m)$.
Similarly to (\ref{4-15.4}) the space
$\overline {\mathcal Z}(N):= \cup \overline {\mathcal Z}_{\leq w}(N) $
is  an algebra \emph{bifiltred} by the weight and by the depth.
We want to describe this algebra and its
associate graded for the weight and depth filtrations.
Let us start from some relations between these numbers.
Notice that
$$
Li_1(\zeta_N^{\alpha}) = -\log(1-\zeta_N^{\alpha})
$$
so the simplest case of this problem reduces to theory of cyclotomic units. By  the Bass theorem
all relations between the cyclotomic units $1-\zeta_N^{\alpha}$ follow from the distribution
relations and symmetry under $\alpha \to -\alpha$ valid modulo roots of unity.

\subsection{Relations}

\emph{The double shuffle relations}.
Consider the generating series
$$
Li (x_{1},\dots ,x_{m}|t_1,\dots ,t_m ) :=
\sum_{n_i \geq 1}Li_{n_{1},\dots ,n_{m}}(x_{1},\dots ,x_{m} )t_1^{n_1-1}\dotsb  t_m^{n_m-1}\,.
$$
Let $\Sigma_{p,q}$  be the subset of permutations of $p+q$ letters $\{1,\dots ,p+q\}$
consisting of all shuffles of  $\{1,\dots ,p\}$ and $\{p+1,\dots ,p+q\}$.
Similarly to (\ref{4-15.3}) multiplying  power series (\ref{zhe5}) we immediately get
\begin{multline} \label{shuffle1}
    Li(x_1,\dots ,x_{p}|t_1,\dots ,t_{p})\cdot
    Li(x_{p+1},\dots ,x_{p+q}|t_{p+1},\dots ,t_{p+q}) =\\
    = \sum_{ \sigma \in \Sigma_{p,q}}Li(x_{\sigma (1)},\dots ,x_{\sigma (p+q) }|t_{\sigma ( 1) },
    \dots , t_{\sigma (p+q)}) \quad + \quad \mbox{lower depth terms}\,.
\end{multline}
To get the other set of the relations we  multiply iterated integrals (\ref{2/7/00/1}),
 and use theorem~\ref{BB} plus the following product formula
for the iterated integrals:
\begin{equation} \label{4-15.5}
\int_{\gamma} \omega_{1} \circ \dots  \circ \omega_{p} \cdot
\int_{\gamma} \omega_{p+1} \circ \dots  \circ \omega_{p+q}  \quad = \quad \sum_{ \sigma \in \Sigma_{p,q}}
\int_{\gamma} \omega_{\sigma(1)} \circ \dots  \circ \omega_{\sigma(p+q)}\,.
\end{equation}
For example the simplest case of formula (\ref{4-15.5}) is derived as follows:
$$
\int_0^1 f_1(t)dt \cdot \int_0^1f_2(t)dt  \quad = \quad \Bigl(\int_{0 \leq t_1 \leq t_2 \leq 1} +
\int_{0 \leq t_2 \leq t_1 \leq 1}\Bigr)f_1(t_1)f_2(t_2)dt_1dt_2\,.
$$
It is very similar in spirit to the derivation (\ref{4-15.4}) of formula (\ref{4-15.3}).

To get nice formulas consider the generating series
\begin{multline} \label{4-15.30}
I^*(a_{1}:\dots :a_{m}:a_{m+1}|t_1,\dots ,t_m ) :=\\
\hspace*{-2mm}\sum_{n_i \geq 1}I_{n_{1},\dots ,n_{m}}(a_{1}:\dots :a_{m}:a_{m+1})
t_1^{n_1-1} (t_1+t_2)^{n_2-1}\!\dotsb(t_1+\dots  +t_m)^{n_m-1}\,.
\end{multline}

\begin{theorem}
\begin{multline} \label{shuffle2}
I^*( a_1:\dots :a_{p}:1 |t_1,\dots ,t_{p})\cdot
I^*( a_{p+1}:\dots :a_{p+q}:1 |t_{k+1},\dots ,t_{p+q}) =\\
=\sum_{\sigma \in \Sigma_{p,q}}I^*(a_{\sigma (1)},\dots ,a_{\sigma (p+q)} :1 |t_{\sigma ( 1) },\dots ,t_{\sigma (p+q)})\,.
\end{multline}
 \end{theorem}

\smallskip\noindent
\textit{A sketch of the proof}
It is not hard to prove the following formula
\begin{equation} \label{43212}
I^*(a_1:\dots :a_m:1|t_1,\dots ,t_m) = \int_0^1\frac{ s^{-t_1}}{a_1 -
  s}ds \circ \dots  \circ \frac{s^{-t_m}}{a_m - s}ds\,.
\end{equation}
The theorem follows   from this and
product formula (\ref{4-15.5}) for the  iterated integrals.

For multiple $\zeta$'s these are precisely
the  relations of Zagier, who conjectured that, properly regularized, they
provide all the relations between the multiple $\zeta$'s.

\emph{Distribution relations}. From the power series expansion we immediately get

\begin{proposition} \label{2.5}
   If $|x_i| < 1$  and $l$ is a positive integer then
\begin{equation} \label{5n.new}
 Li (x_{1},\dots ,x_{m}|t_{1},\dots ,t_{m})=\sum_{ y^l_i
 =  x_i} Li( y_{1},\dots , y_{m}|lt_{1},\dots ,lt_{m})\,.
\end{equation}
\end{proposition}

If $N>1$ the double shuffle plus distribution relations do not
provide all relations between multiple polylogarithms at $N$-th roots of unity.
However I conjecture they do give all the relations
if $N$ is a prime and we restrict to  the weight $=$ depth case.

\subsection{Multiple polylogarithms at roots of unity and the cyclotomic  Lie algebras}

Denote by  $UC_{\bullet}$  the universal enveloping algebra
of a graded Lie algebra $C_{\bullet}$. Let $UC_{\bullet}^{\vee}$ be its graded dual.
It is a commutative Hopf algebra.

\begin{conjecture} \label{cycle}
\begin{enumerate}[a)]
\item There exists a graded Lie algebra $C_{\bullet}(N)$ over $\Q$
   such that one has an isomorphism $
\overline {\mathcal Z}(N) \quad = \quad UC_{\bullet}(N)^{\vee}
$
of filtered  by the weight on the left and by the degree on the right
algebras.

\item $H^1_{(n)}(C_{\bullet}(N)) \quad = \quad
K_{2n-1}(\Z[\zeta_N][\frac{1}{N}]) \otimes \Q$.

\item $C_{\bullet}(N) \otimes\Q_l = {\mathcal G}_N^{(l)}$ as filtered by the weight Lie algebras.
\end{enumerate}
\end{conjecture}

Here $H_{(n)}$ is the degree $n$ part of $H$.
Notice that $H^1_{(n)}(C_{\bullet}(N))$ is dual to the space of
degree $n$ generators of the Lie algebra $C_{\bullet}(N)$.

\begin{examples}
\begin{enumerate}[i)]
\item If $N=1$ the generators should correspond to $\overline\zeta(2n+1)$.
\item If $N>1, n>1$ the generators should correspond
$ \overline Li_{n }( \zeta_N^{\alpha} )$ where $(\alpha, N)=1$.
\end{enumerate}
\end{examples}

A construction of the Lie algebra $C_{\bullet}(N)$ using the Hodge theory see in \cite{G3}. Similarly to theorem 1.2 one  proves (\cite{G8}) that the algebra 
$\overline {\mathcal Z}(N)$ is a subalgebra of the universal enveloping 
algebra of the motivic Tate Lie algebra of the scheme $S_N$, i.e. 
free graded Lie algebra generated by $K_{2n-1}(S_N) \otimes \Q$ in degrees $n \geq 1$. However this estimate is not exact for sufficiently big $N$. 

\subsection{The coproduct}

The iterated integral
\begin{equation} \label{4-30.3}
{I}(a_0; a_1,  \dots , a_m; a_{m+1}):= \int_{a_0}^{a_{m+1}}\frac{dt}{t-a_1}\circ \dots
\circ \frac{dt}{t-a_m}
\end{equation}
provides a framed mixed Hodge-Tate structure,
denote by ${I}_{{\mathcal H}}(a_0; a_1,\! \dots\! , a_m; a_{m+1})$,
 see \cite{G1,G3}.
The set of equivalence classes of framed mixed Hodge-Tate structures has a structure of
the graded
Hopf algebra over $\Q$ with the coproduct $\Delta$.

\begin{theorem} \label{CP1} For the
framed Hodge-Tate structure corresponding to
(\ref{4-30.3}) we have:
\begin{multline} \label{CP2}
\Delta {I}_{{\mathcal H}}(a_0; a_1, a_2, \dots , a_m; a_{m+1})=\\
= \!\!\sum_{0 = i_0 < i_1 < \dots  < i_k < i_{k+1} = m}\hspace{-11mm} {I}_{{\mathcal H}}(a_0; a_{i_1}, \dots ,
a_{i_k}; a_{m+1}) \otimes \prod_{p =0}^k
{I}_{\mathcal H}(a_{i_{p}}; a_{i_{p}+1}, \dots , a_{i_{p+1}-1}; a_{i_{p+1}})
\end{multline}
\end{theorem}

This formula also provides an explicit description
of the variation of mixed Hodge-Tate structures whose period function is
given by (\ref{4-30.3}), see \cite{G3,G8}.
Specializing it we get explicit formulas for the coproduct
of all multiple polylogarithms (\ref{2/7/00/1}).
When $a_i$ are  $N$-th roots of unity and the lower depth terms are suppressed
the result has a particular nice form. It is described, in an axiomatized
form of the coproduct for dihedral Lie algebras, in the next section.

\section{The Dihedral Lie Coalgebra of a Commutative Group \boldmath $G$}

Let $G$ and $H$ be two commutative groups or, better,
commutative group schemes. Then, generalizing a construction given in \cite{G4,G5}
one can define a graded Lie coalgebra ${\mathcal D}_{\bullet}(G|H)$, called the dihedral
Lie coalgebra of $G$ and $H$ (\cite{G8}).
In the special case when $H = {\rm Spec}\Q[[t]]$ is the additive group of the
formal line it is a bigraded Lie coalgebra
${\mathcal D}_{\bullet \bullet}(G)$ called the dihedral Lie
coalgebra of  $G$. (The second grading is coming from the natural
filtration on $\Q[[t]]$). We recall its definition below. The construction of  ${\mathcal D}_{\bullet}(G|H)$ is left as an
 easy exercise.

\subsection{Formal definitions (\cite{G4,G5})}

Let $G$ be a commutative group.
We will define a bigraded Lie coalgebra
${\mathcal D}_{\bullet \bullet}(G ) = \oplus_{w \geq m \geq 1}
 {\mathcal D}_{w, m}(G )$.
The $\Q$-vector space  $ {\mathcal D}_{w,m}(G)$ is generated by the symbols
\begin{equation} \label{ginf}
I_{n_1,\dots ,n_m}(g_1: \dots  : g_{m+1})\,,    \qquad w =
n_1+\dots +n_m, \quad n_i \geq 1\,.
\end{equation}
To define the relations we introduce the generating series
\begin{multline} \label{ginf11}
\{g_1: \dots  : g_{m+1}|t_1:\dots :t_{m+1}\}:=\\
\sum_{n_i >0} I_{n_1,\dots ,n_m}(g_1: \dots  : g_{m+1})
(t_1-t_{m+1})^{n_1-1}\dotsb(t_m -t_{m+1})^{n_m-1}\,.
\end{multline}

We will also need two other generating series:
 \begin{equation} \label{ginf112}
    \{g_1: \dots  :  g_{m+1}|t_1,\dots ,t_{m+1}\} :=
    \{g_1:   \dots  : g_{m+1}
    |t_1:t_1+t_2:\dots :t_1+\dots +t_{m}:0 \}
\end{equation}
where $t_1+ \dots + t_{m+1} =0$, and
\begin{equation} \label{::}
\{g_1,   \dots  ,  g_{m+1}|t_1:\dots :t_{m+1}\} := \{1: g_1: g_1  g_2 : \dots  :
g_1   \dotsb   g_{m }|
  t_1: \dots  :  t_{m+1} \}
\end{equation}
   where $g_1 \cdot  \dotsb  \cdot  g_{m+1} = 1$.

\subsection{Relations}
\begin{enumerate}[i)]
\item \emph{Homogeneity}. For any $g \in G$ one has
\begin{equation} \label{gw0}
\{g\cdot g_1:   \dots  : g\cdot g_{m+1}|t_1: \dots  :t_{m+1}\} =
\{g_1:  \dots  : g_{m+1}|t_1: \dots : t_{m+1}\}\,.
\end{equation}
(Notice that the homogeneity in $t$ is true by the very definition (\ref{ginf11})).

\item \emph{The double shuffle relations}  $(p+q = m, p\geq 1, q \geq 1)$.
\begin{equation} \label{gw3}
 \sum_{\sigma \in \Sigma_{p,q}}\{ g_{\sigma(1)} : \dots  :  g_{\sigma(m)}: g_{m+1}|   t_{\sigma(1)}, \dots ,  t_{\sigma(m)}, t_{m+1}\} =0,
\end{equation}
\begin{equation} \label{gw4}
 \sum_{\sigma \in \Sigma_{p,q}}
\{  g_{\sigma(1)}, \dots  , g_{\sigma(m)}, g_{m+1} | t_{\sigma(1)}: \dots :
t_{\sigma(m)}: t_{m+1} \} =0\,.
\end{equation}

\item \emph{The distribution relations}. Let $l \in \Z$. Suppose that the $l$-torsion subgroup $G_l$ of  $G$ is finite
and its order is divisible by $l$. Then if $x_1, \dots , x_m$ are $l$-powers
\begin{multline*}
\{x_1:  \dots  : x_{ m+1}| t_1:  \dots :  t_{ m+1 } \} -\\
 - \frac{1}{|G_l|}\sum_{y_i^l = x_i}\{y_1: \dots  :y_{m+1}| l\cdot t_1: \dots :  l\cdot t_{ m+1} \} =0
\end{multline*}
except the relation  $I_1(e:e) =   \sum_{y^l=e} I_1(y:e)$ which
  is not supposed to  hold.

\item $I_1(e:e) = 0$.
\end{enumerate}

Denoted by $\widehat {\mathcal D}_{\bullet \bullet}(G )$ the bigraded space defined just as above except condition iv) is dropped, so $\widehat {\mathcal D}_{\bullet \bullet}(G)
= {\mathcal D}_{\bullet \bullet}(G ) \oplus \Q_{(1,1)}$ where $\Q_{(1,1)}$ is of bidgree $(1,1)$.

\begin{theorem} \label{2/7/00/3} {\rm (See theorem 3.1 in
\cite{G5}}.
  If $m \geq 2$ the double shuffle relations imply the
dihedral symmetry relations, which include the cyclic symmetry
$$
\{g_1: g_2: ... :  g_{m+1}|t_1:t_2:...:t_{m+1}\} =
\{g_2: ... : g_{m+1}: g_1|t_2:...:t_{m+1}: t_1\}
$$
the reflection relation
$$
  \{g_1: ... :g_{m+1}|t_1:...:t_{m+1}\} =
(-1)^{m+1}\{g_{m+1}:  ... : g_1|-t_{m}:...:-t_1: -t_{m+1}\}
$$
and the inversion relations
$$
  \{g_1: ... :g_{m+1}|t_1:...:t_{m+1}\} =
\{g_{1}^{-1}:  ... : g_{m+1}^{-1}|-t_{1}:...: -t_{m+1}\}
$$
\end{theorem}

\subsection{Pictures for the definitions}

We think about generating series
(\ref{ginf11}) as a function of $m+1$ pairs
$(g_1,t_1), \dots  , $ \linebreak $ (g_{m+1},t_{m+1})$ located cyclically on an oriented circle as follows.
The oriented circle has slots, where the $g$'s sit, and in between the consecutive slots,
dual slots, where $t$'s sit:

\begin{center}
\ \\
\epsffile{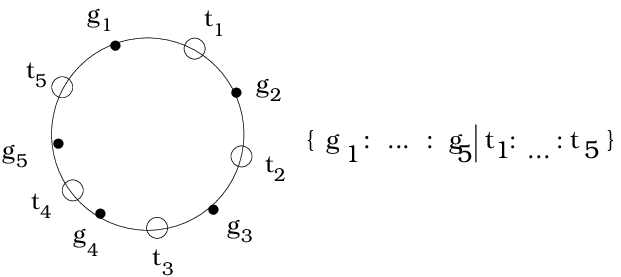}
\ \\
\end{center}

To make definitions (\ref{ginf112}) and (\ref{::})
more transparent set  $g_i' := g_i^{-1}g_{i+1}$, \linebreak $t_i':= -t_{i-1} + t_i$
and put them on the circle together with $g$'s and $t$'s as follows:

\begin{center}
\ \\
\epsffile{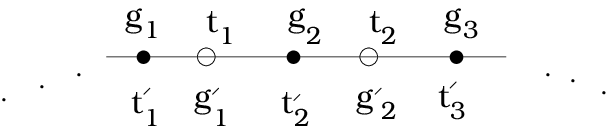}
\ \\
\end{center}
Then it is easy to check that
\begin{multline} \label{9.15.99.1}
\{g_1: \dots  :  g_{m+1}|t_1: \dots  : t_{m+1}\}  =
\{g_1: \dots  :  g_{m+1}|t_1', \dots  , t_{m+1}'\} =\\
=\{g'_1, \dots  ,  g'_{m+1}|t_1: \dots  : t_{m+1}\}\,.
\end{multline}
To  picture any of three generating series (\ref{9.15.99.1})
we leave on the circle only the two sets of variables among
$g$'s, $g'$'s, $t$'s, $t'$'s which appear in this generating series:
\begin{center}
\ \\
\epsffile{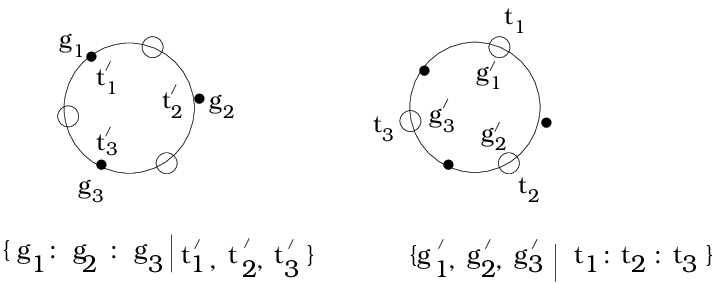}
\ \\
\end{center}
The  ``$\{:\}$''-variables are outside, and the ``$\{,\}$''-variables are inside of the circle.

\subsection{Relation with multiple polylogarithms when \boldmath $G = \mu_N$}

\begin{theorem}
There is a well defined homomorphism of the $\Q$-vector spaces
$$
  {\mathcal D}_{w, m}(\mu_N) \lra {\rm Gr}^{W,D}_{w, m}{\mathcal P}\overline {\mathcal Z}(N)
$$
defined on the generators by
\begin{equation} \label{g21gr}
I_{n_1, \dots , n_m}(a_1: a_2:  \dots  : a_{m+1}) \lra (2\pi i)^{-w}
\mbox{{\rm integral (\ref{2/7/00/1})}}\,.
\end{equation}
\end{theorem}

Here we consider iterated integrals (\ref{2/7/00/1})
modulo the lower depth integrals and products of similar integrals.
Formula (\ref{::}) reflects  theorem~\ref{BB}. Formula
(\ref{ginf112}) reflects definition (\ref{4-15.30}) of the generating series $I^*$.
The shuffle relations (\ref{gw3}) (resp. (\ref{gw4})) correspond to relations (\ref{shuffle1})
(resp. (\ref{shuffle2})).

\begin{remark}
Notice an amazing symmetry between $g$'s and $t$'s in the generating series
(\ref{9.15.99.1}), completely unexpected from the point of view of
iterated integrals (\ref{2/7/00/1}).
\end{remark}

\subsection{The cobracket \boldmath{
$ \delta:     {\mathcal D}_{\bullet \bullet}(G)  \longrightarrow
  \Lambda^2{\mathcal D}_{\bullet \bullet}(G)$}}

It will be defined by
\begin{equation} \label{ccc3}
\begin{split}
&\delta \{g_1:  \dots  : g_{m+1}| t_1: \dots  :t_{m+1}\} =\\
& = -\sum_{k=2}^{m} {\rm Cycle}_{m+1}\left(\{g_{1}:\dots
:g_{k-1}:g_k| t_{1}: \dots  : t_{k-1}: t_{m+1}    \}\right.\\
& \left.\wedge \{g_{k}  :   \dots  : g_{m+1}| t_k: \dots  : t_{m+1}  \} \right)
\end{split}
\end{equation}
where indices are modulo $m+1$ and
${\rm Cycle}_{m+1} f(v_1,\dots ,v_{m}):= \sum_{i=1}^{m+1}f(v_{i},\dots ,$ \linebreak $v_{i+m})$.

Each term of the formula corresponds to the following procedure: choose a slot and
a dual slot on the circle. Cut the circle
 at the chosen  slot and dual slot  and make two   oriented circles
with the data  on each of them obtained from the initial data.
It is useful to think about the slots and dual slots as of little arcs, not points,
so cutting one of them we get the
arcs on each of the two new circles marked by the corresponding letters.
The formula reads as follows:
$$
\delta (\ref{ccc3}) = - \sum_{{\rm cuts}} \mbox{(start at the dual slot)}
 \wedge \mbox{(start at the slot)}
$$

\begin{center}
\ \\
\epsffile{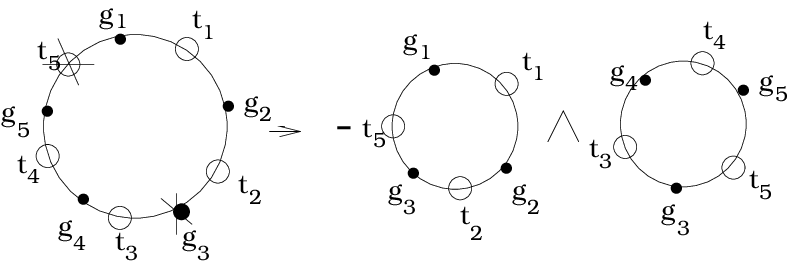}
\ \\
\end{center}

\begin{theorem} \label{9.99.18}
There exists unique map $\delta: {\mathcal D}_{\bullet \bullet}(G) \lra
\Lambda^2 {\mathcal D}_{\bullet \bullet}(G)$ for which (\ref{ccc3}) holds, providing
a bigraded Lie coalgebra structure    on
  ${\mathcal D}_{\bullet \bullet}(G)$.

A similar result is true for $\widehat {\mathcal D}_{\bullet \bullet}(G)$. Moreover there is an isomorphism of bigraded Lie algebras
$
\widehat {\mathcal D}_{\bullet \bullet}(G) =  {\mathcal D}_{\bullet \bullet}(G)\oplus \Q_{(1,1)}
$.
\end{theorem}

\section{The Dihedral Lie Algebra of \boldmath$\mu_N$ and Galois Action on
\boldmath $\pi^{(l)}_{1}(X_N)$ }

\subsection{Constraints on the image of the Galois group}

Recall the   homomorphism
\begin{equation} \label{121212*}
\phi^{(l)}_{N}\colon {\rm Gal}(\overline{\Bbb Q}/ {\Bbb Q}) \longrightarrow
{\rm Aut} {\Bbb L}_N^{(l)}
\end{equation}
It has the following properties.
\begin{enumerate}[i)]
\item \emph{The action of $\mu_N$.} The group $\mu_N$ acts on $X_N$ by $z\lms \zeta_N z$.
This action does not preserve
the base vector $v_{\infty}$. However one can define an action of $\mu_N$ on
$\pi_1^{(l)}(X_N, v_{\infty}) \otimes \Q_l$, 
and hence on ${\Bbb L}^{(l)}(X_N)$,
commuting with the action of the Galois group ${\rm Gal}_{\Q(\zeta_{l^{\infty}N})}$ \cite{G5}.

\item \emph{``Canonical generator'' at $\infty$}. Recall the projection
\begin{equation} \label{PP}
 {\Bbb L}^{(l)}(X_N)\lra
 {\Bbb L}^{(l)}({\Bbb G}_m)   = \Z_l(1)
\end{equation}
Let $X$ be a regular curve over $\overline \Q$,
$\overline X$ the corresponding projective curve and $v$ a tangent vector at
$x \in \overline X$. Then  there is a natural map of Galois modules ([D]):
$$
\Z_l(1) = \pi_1^{(l)}(T_x \overline X \backslash 0, v)  \lra  \pi_1^{(l)}(X, v)
$$
For
$X = {\Bbb G}_m$, $x = \infty, v = v_{\infty}$
it is an isomorphism. So for $X = X_N$ it
provides a
  splitting  of  (\ref{PP}):
\begin{equation} \label{III}
X_{\infty}\colon    \Z_l(1) = {\Bbb L}^{(l)}({\Bbb G}_m)  \hookrightarrow {\Bbb L}^{(l)}(X_N)
\end{equation}

\item \emph{Special equivariant generators for ${\Bbb L}^{(l)}_N$}.
For topological reasons there are well defined
conjugacy classes of ``loops
around $0$ or $\zeta \in \mu_N$'' in $\pi_1(X_N, v_{\infty})$. It turns out
that there are $\mu_N$-equivariant
representatives of these classes   providing a set of the
generators for the Lie algebra ${\Bbb L}^{(l)}_N$:
\end{enumerate}

\begin{lemma}  \label{splitting}
There exist maps $X_0, X_{\zeta}: \Q_l(1) \lra {\Bbb L}^{(l)}_N$
which belong to the conjugacy classes  of the
``loops around $0, \zeta$'' such that
$
X_0 + \sum_{\zeta \in \mu_N} X_{\zeta} + X_{\infty} =0
$ and the action of $\mu_N$ permutes $X_{\zeta}$'s (i.e.\ $\xi_*X_{\zeta} = X_{\xi \zeta}$)
and fixes $X_0$, $X_{\infty}$.
\end{lemma}

To incorporate these constraints we employ a more general set up.

\subsection{The Galois action and special equivariant derivations}

Let $G$ be a commutative group written multiplicatively.
Let  ${L}(G)$ be the free Lie algebra with the generators $X_i$ where $i \in \{0\} \cup G$
(we assume  $0 \not \in G$). Set $X_{\infty}:= -X_0 -\sum_{g \in G}X_g$.

A derivation ${D}$ of the Lie algebra
${L}(G)$ is called special if there are elements  $S_i \in {L}(G)$ such that
\begin{equation} \label{g2}
{D}(X_{i}) = [S_{i}, X_{i}] \quad \mbox{for any $i \in \{0\} \cup G$},
\quad \mbox{and} \quad {D}(X_{\infty}) = 0\,.
\end{equation}
The special derivations of $L(G)$ form a Lie algebra,  denoted ${\rm Der}^{S}{L}(G)$.
Indeed, if $D(X_i) = [S_i,X_i]$, $D'(X_i) = [S'_i,X_i]$, then
\begin{equation} \label{DDD}
[D,D'](X_i) = [S''_i, X_i], \quad \mbox{where} \quad S''_i :=
D(S'_i) - D'(S_i) +[S'_i, S_i]\,.
\end{equation}
The group $G$ acts on the generators by $h\colon X_0 \lms X_0$, $X_g \lms X_{hg}$.
So it acts by  automorphisms of the Lie algebra ${L}(G)$. A derivation ${D}$ of
${L}(G)$ is called equivariant if it  commutes with
the action of $G$.
 Let ${\rm Der}^{SE}{L}(G)$ be the Lie algebra of all special equivariant derivations
of the Lie algebra ${L}(G)$.

Denote by ${\Bbb L}(\mu_N)$ the pronilpotent completion
of the Lie algebra ${L}(\mu_N)$.
Lemma~\ref{splitting} provides a (non canonical) isomorphism
${\Bbb L}(\mu_N) \otimes \Q_l \lra {\Bbb L}^{(l)}_N$. Then it follows
that ${\mathcal G}_N^{(l)}$ acts by special equivariant derivations of the Lie algebra
${\Bbb L}_N^{(l)}$, i.e.\
\begin{equation} \label{MYF}
{\mathcal G}_N^{(l)} \hookrightarrow  {\rm Der}^{SE}{\Bbb L}^{(l)}_N\,.
\end{equation}

\subsection{Incorporating the two filtrations}

The Lie algebra $L(G)$ is bigraded
by the weight and depth. Namely, the free generators $X_0$, $X_{g}$ are bihomogeneous:
they are of weight $-1$, $X_0$ is of depth $0$ and the $X_g$'s are of depth $-1$.
Each of the gradings induces a filtration of $L(G)$.

The Lie algebra ${\rm Der}{L}(G)$ is bigraded by the weight and depth. Its Lie subalgebras
 ${\rm Der}^S{L}(G)$ and ${\rm Der}^{SE}{L}(G)$ are compatible with the
weight grading. However they are \emph{not} compatible with the depth
grading. Therefore they are graded by the weight, and \emph{filtered} by
the depth. A derivation (\ref{g2}) is of depth $-m$ if
 each $S_j$ mod $X_j$ is of depth $-m$, i.e.\ there are at least $m$ $X_i$'s different
from $X_0$ in $S_j$ mod $X_j$. The depth filtration is compatible with the weight grading.
Let ${\rm Gr}{\rm Der}^{SE}_{\bullet \bullet}{L}(G)$ be the
associated graded for the depth filtration.

Consider the following linear algebra situation.
Let $W_{\bullet}L$ be a filtration on a vector space $L$. A splitting $\varphi\colon
{\rm Gr}^W L \lra L$ of the filtration leads to an isomorphism
$\varphi^*\colon {\rm End}(L) \lra {\rm End}({\rm Gr}^WL)$. The space ${\rm End}(L)$
inherits a natural filtration, while ${\rm End}({\rm Gr}^WL)$ is graded.
The map $\varphi^*$ respects the corresponding filtrations.
The map
$ {\rm Gr}\varphi^*\colon {\rm Gr}^W({\rm End}L)  \lra
{\rm End}({\rm Gr}^WL)
$ does not depend on the choice of the splitting.
Therefore if  $L = {\Bbb L}_N$ we get a \emph{canonical} isomorphism
\begin{equation} \label{PPPQ}
{\rm Gr}^W({\rm Der}^{SE}{\Bbb L}_N) \stackrel{\sim}{=}{\rm Der}^{SE}{L}(\mu_N)
\end{equation}
respecting the weight grading.
Thus there is a \emph{canonical}
injective morphism
\begin{equation} \label{myLI}
{\rm Gr}{\mathcal G}^{(l)}_{\bullet \bullet}(\mu_N)  \stackrel{}{\hookrightarrow}
{\rm Gr}{\rm Der}^{SE}_{\bullet \bullet}{\Bbb L}^{(l)}_N
 \stackrel{\sim}{=} {\rm Gr}{\rm Der}^{SE}_{\bullet \bullet}{L}(\mu_N) \otimes \Q_l\,.
\end{equation}

The Lie algebra ${\mathcal G}^{(l)}_N$ is isomorphic to
${\rm Gr}_{\bullet}^W{\mathcal G}^{(l)}_N$,
but this isomorphism is not canonical. The advantage of working
with ${\rm Gr}_{\bullet}^W{\mathcal G}^{(l)}_N$ is that, via  isomorphism (\ref{PPPQ}),
 it became a Lie subalgebra of
${\rm Der}_{}^{SE}{L}(\mu_N) \otimes \Q_l$, which has natural generators
provided by the canonical generators of ${L}(\mu_N)$. This gives canonical
``coordinates'' for description of  ${\rm Gr}_{\bullet}^W{\mathcal G}^{(l)}_N$.
The  benefit of taking its associated graded for the \emph{depth} filtration is an
unexpected relation with the geometry of modular varieties for $GL_m$, where $m$ is the depth.

\subsection{Cyclic words and special differentiations (\cite{Dr,K})}

Denote by $A(G)$ the free associative algebra generated by elements $X_i$ where
$i \in \{0\} \cup G$.
Let ${\mathcal C}(A(G)):= A(G)/[A(G), A(G)]$ be the space of cyclic words in $X_i$.
Consider a map of linear spaces
$\partial_{X_i}: {\mathcal C}(A(G)) \lra A(G)$
given on the generators by the following formula (the indices are modulo $m$):
$$
\partial_{X_j} {\mathcal C}(X_{i_1} \dotsb  X_{i_m}):= \sum_{X_{i_k} = X_j} X_{i_{k+1}} \dotsb
X_{i_{k+m-1}}\,.
$$
For example $\partial_{X_1}{\mathcal C}(X_{1} X_{2}X_{1} X^2_{2}) = X_2 X_1 X_2^2 + X_2^2 X_1 X_2$.

Define special derivations just as in (\ref{g2}), but with $S_i \in A(G)$.
There is a map
$$
\kappa: {\mathcal C}(A(G)) \lra {\rm Der}^SA(G), \quad \kappa{\mathcal C}(X_{i_1} \dots  X_{i_m})(X_j):=
[\partial_{X_j} {\mathcal C}(X_{i_1}, \dots , X_{i_m}), X_j]\,.
$$
It is easy to check that it is indeed a special derivation. Denote by
$\widetilde {\mathcal C}(A(G))$ the quotient of ${\mathcal C}(A(G))$ by the subspace generated by the
monomials $X_i^n$. Then one can show that the map $\kappa$ provides an isomorphism of
vector spaces
$$
\kappa\colon \widetilde {\mathcal C}(A(G)) \lra {\rm Der}^SA(G)\,.
$$

\subsection{The dihedral Lie algebra as a Lie subalgebra of special equivariant derivations}

Let $G$ be a finite commutative group. We will use a notation $Y$ for the generator $X_0$ of $A(G)$.
So $\{Y, X_g\}$ are the generators of the algebra $A(G)$. Set
$$
{\mathcal C}(X_{g_0} Y^{n_0-1} \cdot \dotsb  \cdot
X_{g_m} Y^{n_m-1})^G:= \sum_{h \in G} {\mathcal C}(X_{hg_0} Y^{n_0-1} \cdot \dotsb  \cdot
X_{hg_m} Y^{n_m-1})\,.
$$
Consider the following formal expression:
\begin{equation} \label{cw2}
\xi_G :=\sum \frac{1}{|{\rm Aut}{\mathcal C}|}
    I_{n_0, \dots , n_m}(g_0: \dots  : g_m) \otimes {\mathcal C}(X_{g_0} Y^{n_0-1} \cdot \dotsb  \cdot
X_{g_m} Y^{n_m-1})^G
\end{equation}
where the sum is over all  $G$-orbits on the set of cyclic words ${\mathcal C}$ in $X_g , Y$.
The weight $1/|{\rm Aut}{\mathcal C}|$ is the order of automorphism group of the
cyclic word ${\mathcal C}$.

Applying the map $\id \otimes  {\rm Gr}({\kappa})$ we get
a  bidegree $(0,0)$ element
$$
\xi_G \in {\mathcal D}_{\bullet \bullet}(G) \widehat \otimes_{\Q}
{\rm Gr}{\rm Der}^{SE}_{\bullet \bullet}{A}(G)\,.
$$
Let ${D}_{-w,-m}(G)= {\mathcal D}_{w,m}(G)^{\vee}$. Then $
{D}_{\bullet \bullet}(G):=
\oplus_{w,m \geq 1} {D}_{-w,-m}(G) $
 is a bigraded Lie algebra.
Consider $\xi_G$ as a map of bigraded spaces:
 \begin{equation} \label{g3}
\xi_G \in Hom_{\Q-Vect}(  {D}_{\bullet \bullet}(G) ,
{\rm Gr}{\rm Der}^{SE}_{\bullet \bullet}{A}(G))\,.
\end{equation}
Notice that ${\rm Gr}{\rm Der}^{SE}_{\bullet \bullet}{L}(G)$ is a Lie subalgebra of
${\rm Gr}{\rm Der}^{SE}_{\bullet \bullet}{A}(G)$.

\begin{theorem} \label{ga}
 The map $\xi_G$ provides an injective Lie algebra morphism
\begin{equation} \label{12.22.99}
\xi_G\colon {D}_{\bullet \bullet}(G) {\hookrightarrow} {\rm Gr}{\rm Der}^{SE}_{\bullet \bullet}
{L}(G)\,.
\end{equation}
\end{theorem}

\begin{theorem} \label{ramie}
$ {\rm Gr}{\mathcal G}_{\bullet \bullet}^{(l)}(\mu_N)  \hookrightarrow
\xi_{\mu_N}({D}_{\bullet \bullet}(\mu_N)) \otimes_{\Q} {\Q_l}
$.
\end{theorem}

If $G$ is a trivial group we set ${D}_{\bullet \bullet} : = {D}_{\bullet \bullet}(\{e\})$,
 and $\xi:= \xi_{\{e\}}$. Denote by ${D}_{\bullet}(G)$ the quotient of
 ${D}_{\bullet \bullet}(G)$ by the components with $w \not = m$.
\begin{conjecture} \label{ramier}
a) One has $\xi({D}_{\bullet \bullet})\otimes \Q_l  = {\rm Gr} {\mathcal G}^{(l)}_{\bullet \bullet} $.

b) Let $p$ be a prime number. Then $
\xi_{\mu_p}({D}_{\bullet}(\mu_p)) =
{\rm Gr}{\mathcal G}^{(l)}_{\bullet}(\mu_p)
$.
\end{conjecture}

Summarizing we see the following picture: both  Lie algebras ${D}_{\bullet \bullet}(\mu_N)$ and
 ${\rm Gr} {\mathcal G}^{(l)}_{\bullet \bullet}(\mu_N)$  are realized as Lie subalgebras
of the Lie algebra of special equivariant derivations ${\rm Gr}{\rm Der}^{SE}_{\bullet \bullet}{L}(G)$.
The (image of) dihedral Lie algebra contains the
(image of) Galois. Hypothetically they coincide
 when $N=1$, or when weight=depth and $N$ is prime. In general the gap between them  exists,
but should not be big.

\begin{theorem} \label{ramier1}
Conjecture~\ref{ramier} is true for $m=1,2,3$.
\end{theorem}

The proof of this theorem is based on the following two ideas:
the standard cochain complex of
${D}_{\bullet \bullet}(\mu_N)$ is related to the  modular complex, and the modular complex
has a geometric realization. We address them in the last two sections.

\section{Modular Complexes and Galois Symmetries of
\boldmath $\pi^{(l)}_{1}(X_N)$}

\subsection{The modular complexes}
Let $L_m$ be a rank $m$ lattice.
The rank $m$ modular complex ${M}^{\bullet}(L_m) =
{M}_{(m)}^{\bullet}$ is a complex of $GL_m(\Z)$-modules
$$
M_{(m)}^{1} \stackrel{\partial}{\lra} M_{(m)}^{2}
\stackrel{\partial}{\lra}
... \stackrel{\partial}{\lra} M_{(m)}^{m}
$$
If $m=2$ it is isomorphic to complex (\ref{dep44}). In general it
is
defined as follows.

i) {\it The group $M_{(m)}^{1}$}.
  An extended basis of a lattice $L_m$ is an  $(m+1)$-tuple of
vectors
$v_1, ..., v_{m+1}$
of the lattice such that  $v_1 + ... + v_{m+1} = 0$ and
$v_1,...,v_m$ is
a basis.
The extended
basis
form  a principal homogeneous space  over $GL_m(\Z)$.

The abelian group $M^{1}(L_m) = M_{(m)}^{1}$ is generated  by
extended
basis. Denote by
$<e_1, ... , e_{m+1}>$ the generator corresponding to the extended
basis $e_1, ... , e_{m+1}$.

Let $u_1, ..., u_{m+1}$ be elements of the lattice $L_m$ such that
the
set of elements
$\{(u_i, 1)\}$ form a basis of  $L_m \oplus \Z$. The lattice $L_m$
acts
on
such sets by
$
l: \{(u_i, 1)\} \lms \{(u_i+l, 1)\}
$.
We call the coinvariants of this action
{\it homogeneous affine basis} of $L_m$ and denote them by $\{u_1:
... :
u_{m+1}\}$.

To list the relations we need another sets of the generators
corresponding to the homogeneous affine
basis of $L_m$ (compare with (\ref{ginf112}):
$$
<u_1: ... :u_{m+1}>:= <u'_1, u'_2, ..., u'_{m+1}>; \qquad u_i':=
u_{i+1}-u_i
$$
We will also employ the notation
$$
[v_1,...,v_k]:= <v_1,...,v_k, v_{k+1}>, \quad v_1 + ... + v_k +
v_{k+1}
= 0
$$
{\bf Relations}.   One has $<v,-v> = <-v,v>$.
For any $1 \leq k \leq m$, $m\geq 2$ one has (compare with
(\ref{gw3}) -
(\ref{gw4})):
\begin{equation} \label{sshh1}
\sum_{\sigma \in \Sigma_{k,m-k}}
<v_{\sigma(1)}, ... ,v_{\sigma(m)}, v_{m+1}> \quad = \quad 0
\end{equation}
  \begin{equation} \label{sshh3}
\sum_{\sigma \in \Sigma_{k,m-k}} <u_{\sigma(1)}: ...
:u_{\sigma(m)}:
u_{m+1}> \quad = \quad 0
\end{equation}

ii) \emph{The group $M_{(m)}^{k}$}. It is the sum of the groups
$
M^{1}(L^1) \wedge \dots  \wedge M^{1}(L^k)
$
over all unordered lattice decompositions $L_m = L^{1} \oplus \dots  \oplus L^{k}$.
Thus it is generated by the elements $[A_1] \wedge \dots  \wedge [A_k]$ where $A_i$ is a basis
of the sublattice $L_i$ and $[A_i]$'s anticommute.
 Define a map $\partial\colon M_{(m)}^{1} \lra M_{(m)}^{2}$
by setting (compare with (\ref{ccc3}))
$$
\partial\colon <v_1,\dots ,v_{m+1}> \lms  -{\rm Cycle}_{m+1}\Bigl(\sum_{k=1}^{m-1}
[v_1,\dots ,v_k] \wedge [v_{k+1},\dots ,v_m] \Bigr)
$$
where  indices are modulo $m+1$. We get the differential in  $M_{(m)}^{\bullet}$
by Leibniz' rule:
$$
\partial( [A_1] \wedge [A_2] \wedge \dotsc ) :=
\partial([A_1]) \wedge [A_2] \wedge \dotsb    -   [A_1] \wedge \partial([A_2]) \wedge \dots
 + \dotsb\,.
$$

\subsection{The modular complexes and the cochain complex of
\boldmath{${D}_{\bullet \bullet}(\mu_N)$}}

Denote by
$\Lambda_{(m, w)}^*{\mathcal D}_{\bullet \bullet}(\mu_N)$ the depth $m$, weight $w$ part of the
standard cochain complex of the Lie algebra
${D}_{\bullet \bullet}(\mu_N)$.

\begin{theorem} \label{dwhn1}
\begin{enumerate}[a)]
\item For  $m>1$ there exists canonical surjective map of  complexes
\begin{equation} \label{dep5}
    \mu^*_{m;w}\colon M^*_{(m)} \otimes_{\Gamma_1(m;N)} S^{w-m}{\rm V}_m  \lra
    \Lambda_{(m, w)}^*{\mathcal D}_{\bullet \bullet}(\mu_N)\,.
\end{equation}

\item Let $N=1$, or $N=p$ is a prime and $w=m$. Then this map is an isomorphism.
\end{enumerate}
\end{theorem}

The map (\ref{dep5}) was defined in  \cite{G4}. Here is the definition when
$w=m$. Notice that
$$
M^*_{(m)} \otimes_{\Gamma_1(m;N)} \Q = M^*_{(m)} \otimes_{GL_m(\Z)}
\Z[\Gamma_1(m;N)\backslash GL_m(\Z)]\,.
$$
The set $\Gamma_1(m;N)\backslash GL_m(\Z)$ is identified with the set
$\{(\alpha_1, \dots , \alpha_m)\}$ of all nonzero vectors in the vector space over $\Z/p\Z$.
Then
$$
\mu^1_{m;m}: [v_1, \dots , v_m] \otimes (\alpha_1, \dots , \alpha_m) \lra
I_{1, \dots , 1}(\zeta_N^{\alpha_1}, \dots ,\zeta_N^{\alpha_{m}} )\,.
$$
The other components $\mu^*_{m;m}$ are the wedge products of the maps $\mu^1_{k;k}$.

\subsection{Modular complexes and cochain complexes of  Galois Lie algebras}

Combining theorems~\ref{dwhn1}a)
 and~\ref{ramie}
 we get a surjective map of complexes
\begin{equation} \label{dep5**}
M^*_{(m)} \otimes_{\Gamma_1(m;N)} S^{w-m}{\rm V}_m  \lra
\Lambda_{(m, w)}^*{\rm Gr}{\mathcal G}^{(l)}_{\bullet \bullet}(\mu_N)^{\vee}\,.
\end{equation}
So using theorem~\ref{dwhn1}b) we reformulate conjecture~\ref{ramier} as follows:

\begin{conjecture} \label{dwhn2} Let $N=1$, or $N=p$  is a prime and $w=m$.
Then the map (\ref{dep5**}) is an isomorphism.
\end{conjecture}

According to theorem~\ref{ramie} the cochain complex of
${D}_{\bullet \bullet}(\mu_N)$ projects onto the cochain complex of the level $N$ Galois Lie algebra.
Combining this with theorem~\ref{ramier1} (and conjecture~\ref{dwhn2}) we  describe the structure of
the Galois Lie algebras
via the modular complexes. Since the modular complexes are
defined very explicitly this leads to a precise description of Galois Lie algebras.
However to get the most interesting results about them
we need the geometric realization of modular complexes, which allows to express the structure
of the Galois Lie algebra in terms of  the \emph{geometry and topology} of modular varieties.

\section{Geometric Realization of Modular Complexes in Symmetric Spaces}

As was emphasized before the rank $m$ modular complex is purely combinatorial object.
Surprisingly it has a \emph{canonical} realization in the symmetric space ${\Bbb H}_n$.
In the simplest case $m=2$
it identifies the rank two modular complex with the chain complex of the modular triangulation of
the hyperbolic plane.

\subsection{Voronoi's cell decomposition of the symmetric space for \boldmath{$GL_m(\R)$}}

Let
$$
{\Bbb H}_m  := GL_m(\R)/O(m)\cdot \R^* \quad = \frac{
 \mbox{$>0$ definite quadratic forms on $V_m^*$}}{\R^*_+}
$$
Let $L_m \subset V_m$ be a lattice.
Any vector $v \in V_m$ defines a nonnegative definite quadratic
form $\varphi(v):= <v,\cdot>^2$ on $V_m^*$.
The convex hull of  the forms $\varphi(l)$,
when $l$ runs through all non zero primitive
vectors of the lattice $L_m$, is
an infinite polyhedra.
Its projection into the closure of
${\Bbb H}_m$ defines a polyhedral decomposition of ${\Bbb H}_m$ invariant under
the symmetry group of the lattice $L_m$.

\begin{example}
If $m=2$ we get the modular triangulation of the
hyperbolic plane.
\end{example}

Denote by $({V}^{(m)}_{\bullet}, d) $ the chain complex of the Voronoi decomposition.

\subsection{The relaxed modular complex}

Consider a version
${\widehat {M}}_{(m)}^{\bullet}$ of the modular complex, called the relaxed modular complex,
where the group ${\widehat  M}^{1}_{(m)}$ is defined using the same generators $[v_1,\dots ,v_m]$
satisfying only the first shuffle relations (\ref{sshh1}) and
the dihedral symmetry relations
$$
    <v_2,\dots ,v_{m+1}, v_1>  = <v_1,\dots ,v_{m+1}> = (-1)^{m+1}<v_{m+1},\dots ,v_1>\,.
$$
The other groups are defined in a similar way. The differential is as before.

\subsection{The geometric realization map}

Denote by $\varphi(v_{1}, \dots  , v_{k})$ the
convex hull of the forms $\varphi(v_{1}), \dots , \varphi(v_{k})$ in the space of quadratic
forms. Let  $v_1,\dots ,v_{n_1}$ and $v_{n_1+1},\dots ,v_{n_1+n_2}$ be two sets of
lattice vectors such that the lattices they generate have zero intersection.
Define the join $\ast$ by
$$
\varphi(v_1,\dots ,v_{n_1}) \ast \varphi(v_{n_1+1},\dots ,v_{n_1+n_2}) :=
\varphi(v_1,\dots , v_{n_1+n_2})
$$
and extend it by linearity. Make a homological complex ${\widehat {M}}^{(m)}_{\bullet}$
out of ${\widehat {M}}_{(m)}^{\bullet}$ by
$$
{\widehat {M}}^{(m)}_{i}:= {\widehat {M}}_{(m)}^{2m-1-i}\,.
$$

\begin{theorem} \label{15}
There exists a canonical  morphism of complexes
$$
\widehat \psi_{\bullet}^{(m)}\colon{\widehat {M}}^{(m)}_{\bullet} \lra
{V}^{(m)}_{\bullet} \qquad \mbox{such that}
$$
\begin{equation} \label{hh}
    \widehat
    \psi_{\bullet}^{(m)}\Bigl([A_1] \wedge \dots  \wedge [A_k] \Bigr) :=
    \widehat \psi_{\bullet}^{(m)}([A_1]) \ast \dots  \ast  \widehat
    \psi_{\bullet}^{(m)}([A_k])\,.
\end{equation}
\end{theorem}
In particular
$\widehat \psi_{m-1}^{(m)}\Bigl([v_1] \wedge \dots  \wedge [v_m] \Bigr)
= \varphi(v_1, \dots  , v_m) $.

To define such a morphism $\widehat \psi_{\bullet}^{(m)}$ one needs only to define
$\widehat \psi_{2m-2}^{(m)}([v_1,\dots ,v_m])$ for vectors $v_1,\dots ,v_m$ forming a
basis of the lattice $L_m$ in such a way that the dihedral and the
first shuffle relations go to zero and
$$
    d \widehat \psi_{2m-2}^{(m)}([v_1,\dots ,v_m]) =
     \widehat  \psi_{2m-3}^{(m)}(\partial [v_1,\dots ,v_m])
$$
where the right hand side is
computed by (\ref{hh}) and the formula for $\partial $.

\subsection{Construction of the map \boldmath{$\widehat \psi_{2m-2}^{(m)}$}}

A \emph{plane tree} is a tree without self intersections
located on the plane. The edges of a tree consist of legs (external edges)
and internal edges. Choose  a lattice $L_m$.
A \emph{colored tree} is a plane tree whose legs are in a bijective correspondence with the elements of an affine basis of the lattice $L_m$.
In particular a colored $3$-valent tree has $2m-1$ edges.
We visualize it as follows:
\begin{center}
\ \\
\epsffile{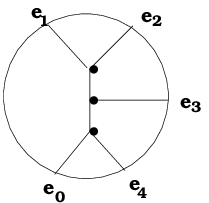}
\ \\
\end{center}
The vectors $e_0,\dots ,e_m$ of an affine basis are located cyclically on an oriented circle
and the legs of the tree end on the circle and labelled by $e_0,\dots ,e_m$.
(The circle itself is not a part of the graph).

\emph{Construction}. Each edge $E$ of the tree $T$ provides a vector $f_E \in L_m$
defined up to a sign. Namely, the edge $E$ determines two trees rooted at $E$, see the picture

\begin{center}
\ \\
\epsffile{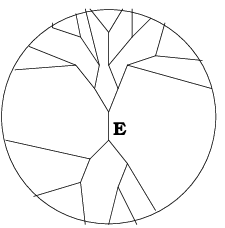}
\ \\
\end{center}

The union of the incoming (i.e.\  different from  $E$) legs
of these rooted trees coincides with the set of all legs of the initial tree.
Take the sum of all the vectors $e_i$ corresponding to the incoming legs  of one of these trees.
Denote it by $f_E$. If we choose the second rooted tree the sum will change the sign.
So the degenerate quadratic form $\varphi(f_E)$ is well defined.
Set
$$
\widehat \psi_{2m-2}^{(m)}(<e_0,e_1,\dots ,e_m>):=
$$
\begin{equation} \label{const.}
\sum_{\mbox{\scriptsize{plane 3-valent trees}}}{\rm sgn} (E_1 \wedge \dots  \wedge E_{2m-1})\cdot
\varphi(f_{E_1}, \dots  , f_{E_{2m-1}})
\end{equation}
Here the sum is over all plane $3$-valent trees colored by $e_0,\dots ,e_m$.
The sign is defined as follows.
Let $V(E)$ be the $\R$-vector space generated by the edges of a  tree. An orientation of a
 tree is a choice of the connected component of ${\rm det}(V(E))\backslash 0$.
A  plane $3$-valent tree has a canonical orientation.
Indeed, the orientation of
the plane provides  orientations of links of each of the vertices.
The sign in (\ref{const.}) is
taken with respect to the canonical orientation of the plane $3$-valent tree.
Then one proves (\cite{G9}) that this map has all the required properties, so we get theorem
\ref{15}.

\begin{examples}
a) For $m=2$ there is one plane $3$-valent tree colored by $e_0, e_1, e_2$, so
we get a modular triangle $\varphi(e_0, e_1, e_2)$ on the hyperbolic plane.
The geometric realization in this case leads to an isomorphism of competes
${M}_{\bullet}^{(2)} \lra {V}_{\bullet}^{(2)} $:
$$
[e_1, e_2] \lms \varphi(e_0, e_1, e_2); \quad [e_1] \wedge [e_2] \lms \varphi(e_1) *  \varphi(e_2) =
\varphi(e_1, e_2)
$$

b) Let $f_{ij}:= e_i+e_j$.
For $m=3$ there are two  plane $3$-valent trees colored by $e_0, e_1, e_2, e_3$,
see the picture, so
the chain is
$$
\widehat \psi_{4}^{(3)}([e_1, e_2, e_3]) := \quad \varphi(e_0, e_1, e_2, e_3, f_{01}) - \varphi(e_0, e_1, e_2, e_3, f_{12})
$$
\begin{center}
\ \\
\epsffile{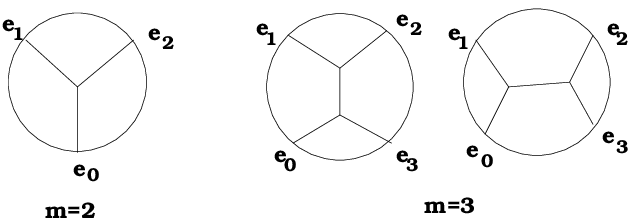}
\ \\
\end{center}
\end{examples}

The symmetric space ${\Bbb H}_3$ has dimension $5$. The Voronoi decomposition consists of the cells of dimensions $5, 4, 3, 2$. All Voronoi cells of dimension $5$ are $GL_3(\Z)$-equivalent to
the Voronoi simplex  $\varphi(e_0, e_1, e_2, e_3, f_{01},  f_{12})$. The map
$\widehat \psi_{4}^{(3)}$ sends the second shuffle relation (\ref{sshh1}) to
the boundary of a Voronoi $5$-simplex.

\begin{theorem} \label{TTCCq}
The geometric realization map provide quasiisomorphisms
$$
{M}_{\bullet}^{(3)} \lra \tau_{[4,2]}({V}_{\bullet}^{(3)});
\qquad {M}_{\bullet}^{(4)} \lra \tau_{[6,3]}({V}_{\bullet}^{(4)})\,.
$$
\end{theorem}

\subsection{Some corollaries}
\begin{enumerate}[a)]
\item Let $N=p$ be a prime.
Take the geometric realization of the rank $3$ relaxed
modular complex. Project it
onto the modular variety
$Y_1(3;p)$. Take the quotient of the group of the $4$-chains generated by
$\widehat \psi^{(3)}_4(v_1, v_2, v_3)$  on $Y_1(3;p)$
by  the subgroup generated by the
boundaries of the Voronoi $5$-cells. Then  the complex
we get is \emph{canonically isomorphic} to the depth=weight $3$ part of the standard cochain complex of
the level $p$ Galois Lie algebra. Therefore
$$
H^i_{(3)}({\rm Gr}\widehat{{\mathcal G}}^{(l)}_{\bullet}(\mu_p))
= H^i(\Gamma_1(3;p)) \quad i=1,2,3\,.
$$
In particular we associate to each of the numbers
$ {Li}_{1,1,1}(\zeta_p^{\alpha_1}, \zeta_p^{\alpha_2}, \zeta_p^{\alpha_3})$,
or to the corresponding Hodge, $l$-adic or  motivic avatars of these numbers,
a certain $4$-cell
on the $5$-dimensional orbifold  $Y_1(3;p)$. The properties of the framed motive
encoded by this number, like the coproduct,  can be read from the geometry of
this $4$-cell.

Similarly the map $\widehat\psi_{2m-2}^{(m)}$ provides a
canonical $(2m-2)$-cell on $Y_1(m;p)$ corresponding to the framed motive with the 
period $Li_{1, \dots  ,1}(\zeta_p^{\alpha_1}, \dots  , \zeta_p^{\alpha_m}) $.

\item $N=1$. Theorems~\ref{ramier1}, \ref{dwhn1}a) and \ref{TTCCq} lead to the following
\end{enumerate}

\begin{theorem} \label{mth1t}
\begin{equation} \label{demfeddd}
{\rm dim} {\rm Gr}{\mathcal G}^{(l)}_{-w, -3} =
\quad    \left\{ \begin{array}{ll}
0 &  w: \quad \mbox{even}  \\
   \left[ \frac{(w-3)^2-1}{48} \right]  &    w:  \quad \mbox{odd}\,. \end{array} \right.
 \end{equation}
\end{theorem}
Since, according to standard conjectures, this number should coincide with $d_{w,3}$
the estimate given in theorem~\ref{1.25.2} should be exact.

\section{Multiple Elliptic Polylogarithms}

The story above is related to the field $\Q$. I hope that for an 
arbitrary number field $F$ there might  be a similar story.  
The Galois group ${\rm Gal}(\overline F/F)$ should have a remarkable 
 quotient ${\rm G}_F$ given by an extension of the maximal abelian quotient of ${\rm Gal}(\overline F/F)$ by a 
prounipotent group $U_F$:
$$
0 \lra U_F \lra {\rm G}_F \lra {\rm Gal}(\overline F/F)^{{\rm ab}}\lra 0
$$ 
Its  structure should be  
related to modular varieties for  $GL_m/_{F}$, for all $m$. 

The group ${\rm G}_{\Q}$ is obtained from 
 the motivic fundamental group of ${\rm G}_m - \{\mbox{``all'' 
roots of unity}\}$. It turns out that for an imaginary quadratic filed $K$ one can get 
a similar picture by taking the motivic  fundamental group of the CM elliptic curve 
$E_K:= \C/{\mathcal O}_K$ punctured at the  torsion points. 
Below we construct the periods of the corresponding mixed motives, {\it 
multiple Hecke L-values}, as the values at torsion points 
of multiple elliptic polylogarithms.  We define the multiple polylogarithms 
for arbitrary curves as 
 correlators for certain Feynman integrals. 
We make sence out of these  Feynman integrals 
by using the perturbation expansion 
via Feynman diagrams,
which in this case are plane $3$-valent trees. Unlike the Feynman integrals 
the coefficients of the perturbative expansion are 
given by convergent finite dimensional integrals, and so well 
defined.   
I leave to the reader the pleasure to penetrate the 
analogy between this construction and the geometric realization 
of modular complexes described in s. 7.4.

\subsection{The classical Eisensten-Kronecker series.} Let $E$ be an elliptic 
curve over $\C$ with the period lattice $\Gamma$, so that $E(\C) = \C/\Gamma$. 
The intersection form $\Lambda^2\Gamma \lra 2 \pi i \Z$ 
leads to the  pairing $\chi: E(\C) \times \Gamma \lra S^1$. 
So for  $a \in E(\C)$ we get a character $\chi_a: \Gamma \lra S^1$. 

Consider the generating function for the classical Eisenstein-Kronecker series
$$
G(a| t):= \quad \frac{ {\rm vol }( \Gamma)}{\pi} \sum'_{\gamma \in \Gamma}
 \frac{\chi_{a}(\gamma)}{|\gamma-t|^2}
$$
where $\sum'$ means the summation over all non zero vectors $\gamma$ of the lattice.  
It depends on a point $a$ of the elliptic curve 
and an element $t$ in  a formal neighborhood of zero in $H_1(E, \R)$. It is invariant
 under the involution $a \lms -a, t \lms -t$. 
Expanding it into power series in $t$ and $\overline t$ we get, as the coefficients, the classical Eisensten-Kronecker series:
$$
G(a| t) = \quad \sum_{p,q \geq 1}\Bigl(\frac{{\rm vol }(\Gamma)}{\pi}  \sum'_{\gamma \in \Gamma}
 \frac{\chi_{a}(\gamma)}{\gamma^p\overline \gamma^q}\Bigr)t^{p-1}\overline t^{q-1}
$$
When $E$ is a CM curve their special values at the torsion points of $E$ 
provide the special values of the Hecke L-series with  
Groessencharacters.

\subsection{Multiple Eisenstein-Kronecker series: a description} 
We define them as 
the coefficients of certain generating functions. The generating function   for 
the depth $m$ multiple Eisenstein-Kronecker series is a 
function
$$
G(a_1: ... :a_{m+1}| t_1, ..., t_{m+1}),\quad t_1 + ... + t_{m+1} =0
$$
where $a_i$ are points on the  elliptic curve $E$ and $t_i$ are elements 
in  a formal neighborhood of zero in $H_1(E, \R)$. 
It is invariant under the shift $a_i \lra a_i +a$. 
Decomposing this function into the series 
in $t_i, \overline t_i$ we get the depth $m$ 
multiple Eisenstein-Kronecker series. 

{\bf Construction}. Consider a plane trivalent tree $T$ 
colored by 
$m+1$ pairs consisting of points $a_i$ on the elliptic curve $E$ 
and formal elements 
$t_i \in H_1(E, \R)$:
\begin{equation} \label{7.23.00.11}
(a_1, t_1), ... , (a_{m+1}, t_{m+1}); \qquad t_1 + ... + t_{m+1}=0
\end{equation}

\begin{center}
\hspace{4.0cm}
\epsffile{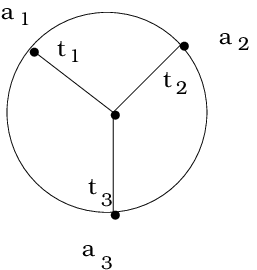}
\end{center}
Each oriented edge $\stackrel{\to}{E}$ of the tree $T$ 
provides an element $t_{\stackrel{\to}{E}} \in H_1(E, \R)$. 
 Namely, as explained in  s. 7.4 the edge $E$ determines two trees rooted at $E$. 
An orientation of the edge $E$ corresponds to the choice of one of them: 
take the tree obtained by going in the direction of 
the orientation of the edge. 
Then $t_{\stackrel{\to}{E}}$ is the sum of all $t_i$'s corresponding to the 
legs  of this tree different from $E$. 
Since $\sum t_i=0$ changing the  orientation of the edge $E$ we 
get  $-t_{\stackrel{\to}{E}}$.

Let $X$ be a manifold and ${\mathcal A}^i(X)$ the space of  $i$-forms on $X$. 
We define  a map
$$
\omega_m: \Lambda^{m+1}{\mathcal A}^0(X) \lra {\mathcal A}^m(X), 
\qquad \omega_m: \varphi_1 \wedge ... \wedge \varphi_{m+1} \lms 
$$
$$
\frac{1}{(m+1)!}{\rm Alt}_{m+1}
\Bigl( \sum_{j=0}^{m+1}(-1)^j \varphi_1 \partial \varphi_2 \wedge ... \wedge \partial \varphi_k \wedge 
\overline \partial \varphi_{k+1} \wedge ... \wedge \overline \partial \varphi_{m+1}\Bigr)
$$
If $ \varphi_i = \log|f_i|$ it is the form used to  define the Chow 
polylogarithms in [G10]. 

Every edge $E$ of the tree $T$ defines a function $G_E$ on $E(\C)^{\{\mbox{vertices of T}\}}$ depending on $t_i$. Namely, let $v_1^E, v_2^E$ be the vertices of the edge $E$. 
Their order orients the edge $E$.  Consider the natural projection 
\begin{equation} \label{7.1.00.22}
p_E: 
E(\C)^{\{\mbox{vertices of the tree T}\}} \quad \lra \quad  E(\C)^{\{v_1^E, v_2^E\}} = E(\C) \times E(\C)
\end{equation}
Then
$
G_E:= \quad p_E^*G \left(x_1^E - x_2^E| t_{\stackrel{\to}{E}}\right)
$ where $(x_1^E, x_2^E)$ is a point at the right of (\ref{7.1.00.22}). 
This function does not depend on the orientation of the edge $E$. 

\begin{definition} 
\begin{equation} \label{7.1.00.2}
G(a_1: ...: a_{m+1}|t_1, ..., t_{m+1}):= 
\end{equation}
$$
 \sum_{\mbox{plane 3-valent trees T}}
 {\rm sgn}(E_1 \wedge ... \wedge E_{2m-1})\cdot\int_{S^{m-1}E(\C)}
{\rm sym}^*\omega_{2m-2}\left(G_{E_1}\wedge ... \wedge 
G_{E_{2m-1}}\right) 
$$
\end{definition}
\noindent
The sum is over all plane $3$-valent trees whose legs are cyclically 
labelled by (\ref{7.23.00.11}). The  correspondence 
$
{\rm sym}: S^{m-1}E(\C)   \lra E(\C)^{\{{\mbox{internal vertices of T}}\}}
$ is given by the  sum  of all $(m-1)!$  natural  maps  $E(\C)^{m-1}\lra E(\C)^{\{{\mbox{internal vertices of T}}\}}$.  

Recall the  CM elliptic curve $E_K$. Let    ${\mathcal N}$ be 
an ideal of ${\rm End}(E_K)$. Denote by $K_{{\mathcal N}}$ the field generated by the 
${\mathcal N}$-torsion points of $E_K$. 
If $a_i$ are ${\mathcal N}$-torsion points of $E_K$
we  view the numbers  obtained in the $t, \overline t$-expansion  (\ref{7.1.00.2}) as multiple 
Hecke L-values related to $K$. They 
 are periods of mixed motives over the ring of integers in 
$K_{{\mathcal N}}$, with ${\rm Norm}({\mathcal N})$ inverted. These are the motives 
which appear in the motivic fundamental group of 
$E - \{{\mathcal N}-\mbox{torsion points}\}$. 

\section{Multiple polylogarithms on curves, Feynman integrals and special values of L-functions}

\subsection{Polylogarithms on curves and special values of L-functions} Let $X$ be a 
regular complex projective algebraic curve  
of genus $g \geq 1$.  Choose a volume form on $X(\C)$, and let $G(x,y)$ be the corresponding Green function. Set 
$$
{\mathcal H}:= H_1(X, \R); \qquad {\mathcal H}_{\C}:= {\mathcal H}\otimes \C = {\mathcal H}^{-1,0} \oplus {\mathcal H}^{0,-1} 
$$
For each integer 
$n \geq 1$ we define a $0$-current $G_n(x,y)$ on $X \times X$ 
with values in 
\begin{equation} \label{7.6.00.1}
{\rm Sym}^{n-1}{\mathcal H}_{\C} (1) \quad = \quad \oplus_{a+b=n-2}
S^{a-1}{\mathcal H}^{-1,0} \otimes S^{b-1}{\mathcal H}^{0,-1} 
\end{equation}
Then $G_1(x,y): = G(x,y)$. For
 $n>1$ it is a function on $X(\C) \times X(\C)$.

To define the 
function $G_n(x,y)$ we proceed as follows. 
Let 
$$
\overline \Omega_a = 
\omega_{\alpha_1}\cdot ... \cdot \omega_{\alpha_{a-1}}\in S^{a-1}\overline \Omega^1, 
\quad 
\Omega_b = \omega_{\beta_1}\cdot ... \cdot \omega_{\beta_{b-1}}\in S^{b-1}\Omega^1; \qquad \omega_{*} \in \Omega^1
$$
Then $\overline \Omega_a \otimes \Omega_b$ is an element of the dual to 
(\ref{7.6.00.1}). We are going to define the  pairing 
$<G_n(x,y), \overline \Omega_a\otimes  \Omega_b >$.  
Denote by $p_i: X^{n-1} \lra X$ the projection on $i$-th factor.

\begin{definition} \label{7.7.00.1} The $n$-th polylogarithm function 
on the curve $X$ is defined by
$$
<G_n(x,y), \overline \Omega_a\otimes \Omega_b > :=\qquad 
$$
$${\rm Alt}_{\{z_1, ..., z_{n-1}\}}\Bigl(\int_{X^{n-1}(\C)}\omega_{n-1}\Bigl(G(x, z_1) \wedge G(z_1, z_2) \wedge ... \wedge G(z_{n-1}, y)\Bigr)  \wedge 
$$
$$
\Lambda_{i=1}^{a-1}p_i^*\overline\omega_{\alpha_i} 
\wedge\Lambda_{j=1}^{b-1}p^*_{a-1+j} \omega_{\beta_j}\Bigr)
$$
\end{definition}
\noindent
We skewsymmetrized the integrand with respect to 
$z_1, ..., z_{n-1}$. 

These functions 
provide a variation of $\R$-mixed Hodge structures on $X \times X - \Delta$ of motivic origin.  
If $X$ is an elliptic curve it is given by  
  Beilinson-Levin theory of elliptic polylogarithms \cite{BL}. 

In particular for a pair of distinct points $x,y$ on $X$
 we get  an $S^{n-1}{\mathcal H}(1)$-framed mixed motive (see \cite{G10} for the background) denoted $\{x,y\}_{n}$, whose period 
is given by $G_n(x,y)$. Its coproduct $\delta$ 
is given by $\{x,y\}_{n} \lms \{x,y\}_{n-1}\wedge (x-y)$ where $(x-y)$ is 
the point of the Jacobian of $X$ corresponding to the divisor $\{x\} - \{y\}$. 
If $X$ is defined over a number field this leads to a very precise conjecture expressing the special value $L(S^{n-1}H^1(X), n)$ 
via the polylogarithms $G_n(x,y)$ - an analog  of  Zagier's conjecture. 
If $X$ is an elliptic curve we are in the situation considered in \cite{G10}, \cite{W}. 
An especially interesting example appears when $x,y$ are cusps on a modular curve. 
Then $(x-y)$ is a torsion point in the Jacobian, so 
$\delta\{x,y\}_n=0$ and thus $G_n(x,y)$ is the regulator of an element of motivic 
$Ext^1(\Q(0), S^{n-1}{\mathcal H}(1))$!

\subsection{Multiple polylogarithms on curves}

We package the polylogarithms $G_n(x,y)$ into the  generating series 
$$
G(x,y| t_1,t_2), \quad t_i \in {\mathcal H}, t_1+t_2=0
$$
so that $G_n(x,y)$ emerges as the weight $-n-1$ component of the 
power series decomposition into $t_1, \overline t_1$. Then  $G(x,y| t_1,t_2) = G(y, x| 
t_2, t_1)$. The construction of the previous section provides 
multiple polylogarithms $G(a_1, ..., a_{m+1}|t_1, ..., t_{m+1})$ on $X$, where 
 $t_1+ ...+ t_{m+1} = 0$. Indeed, for an edge $E$ of a plane $3$-valent 
tree $T$ set $G_E:= p_E^*G(x_1^E, x_2^E| t_{\stackrel{\to}{E}}, - t_{\stackrel{\to}{E}})$ and repeat the construction. We call the constant term in $t$'s the 
multiple Green function on $X$.

\subsection{Feynman integral for multiple Green functions} 
Let $\varphi$ be a function and $\psi$ a $(1,0)$-form 
on $X(\C)$ with values in $N\times N$ complex matrices. We denote by 
$\overline \varphi$ and $\overline \psi$ the result of the action of complex conjugation. Then the multiple Green function $G(a_1, .., a_{m+1})$ emerges 
as the leading term of the asymptotic when $N \to \infty$ of the following correlator:
$$
\int {\rm Tr}\Bigl((\varphi + \overline \varphi)(a_1) \cdot ... \cdot 
(\varphi + \overline \varphi)(a_{m+1})\Bigr) e^{iS(\varphi, \psi)}{\mathcal D}\varphi
{\mathcal D}\psi
$$
where 
$$
S(\varphi, \psi):= \int_{X(\C)} {\rm Tr}
\Bigl(\varphi \overline \partial \psi + \overline \varphi \partial \overline \psi + \psi \overline \psi + \varphi [\psi, \overline \psi ]  + \overline \varphi [\overline \psi, 
\psi ]\Bigr)
$$
 
I conjecture that the special values $L(S^{n}H^1(X), n+m)$ can be expressed via the depth $m$ multiple polylogarithms on $X$. 
    So  Feynman integrals provide construction of (periods of) 
mixed motives, which are in particular responsible for special values of L-functions. 
I hope this reflects a very general phenomena.

\end{document}